\documentclass[journal]{IEEEtran}

\usepackage{booktabs}
\usepackage{float}
\usepackage[colorlinks,linkcolor=blue]{hyperref}
\usepackage[latin1]{inputenc}
\usepackage{amsmath}
\usepackage{graphicx}
\usepackage{amssymb}
\usepackage{amsthm}
\usepackage{verbatim}
\usepackage{epsfig}
\usepackage{enumerate}
\usepackage{blindtext}
\usepackage{epstopdf}
\usepackage{subcaption}
\newcommand{\mrm}{\mathrm}
\newtheorem{theorem}{Theorem}

\newtheorem{assumption}{Assumption}
\newtheorem{lemma}{Lemma}
\newtheorem{problem}{Problem}
\newtheorem{remark}{Remark}

\newtheorem{corollary}{Corollary}
\newtheorem{proposition}{Proposition}
\usepackage[makeroom]{cancel}
\usepackage{stfloats}
\usepackage{amssymb}
\usepackage{amsthm}
\usepackage{graphicx}
\usepackage{float}
\usepackage{amsmath}
\usepackage{amsmath,bm}
\usepackage{color,xcolor}
\usepackage{multirow}
\usepackage{mathrsfs}
\usepackage{lineno,hyperref}
\usepackage{lipsum}
\usepackage{stfloats}
\usepackage{titletoc}
\usepackage{dsfont}
\usepackage{mathtools}

\usepackage{enumitem}
\usepackage{algorithm2e}
\RestyleAlgo{ruled}

\usepackage[prependcaption,colorinlistoftodos]{todonotes}

\newcommand{\nlm}[1]{{\color{blue}{#1}}}

\usepackage[nolist]{acronym}
\begin{acronym}
	\acro{lre}[LRE]{linear regression equation}
	\acro{ltv}[LTV]{linear time-varying}
	\acro{lmi}[LMI]{linear matrix inequality}
	\acro{iss}[ISS]{input-to-state stability}
	\acro{c+i}[C+I]{consensus + innovations}
	\acro{lf}[LF]{Lyapunov function}
	\acro{guas}[GUAS]{Global Uniform Asymptotic Stability}
	\acro{uco}[UCO]{Uniform Complete Observability}
	\acro{pe}[PE]{Persistency of Excitation}
	\acro{cpe}[cPE]{cooperative Persistency of Excitation}
	\acro{sdp}[SDP]{semi-definite program}
\end{acronym}

\DeclareMathOperator{\diag}{diag}
\DeclareMathOperator{\col}{col}

\newcommand{\cPEu}{\bar{\iota}_1}
\newcommand{\cPEl}{\underline{\iota}_1}
\newcommand{\MOIu}{\bar{\iota}_2}
\newcommand{\MOIl}{\underline{\iota}_2}
\newcommand{\Mu}{\bar{\iota}_3}
\newcommand{\Ml}{\underline{\iota}_3}
\newcommand{\la}{\bar{\Lambda}}
\newcommand{\ytil}{\tilde{y}}
\newcommand{\y}{y}
\newcommand{\yhat}{\hat{y}}
\newcommand{\D}{\Delta}
\newcommand{\Db}[1]{\bar{\Delta}_{#1}}
\newcommand{\Th}{x}
\begin{document}
	
	\title{A robust consensus + innovations-based distributed parameter estimator}
	
	\author{Nicolai Lorenz-Meyer, Juan G. Rueda-Escobedo, Jaime~A.~Moreno, and Johannes Schiffer
		\thanks{N. Lorenz-Meyer is with Brandenburg University of Technology Cottbus-Senftenberg, 03046 Cottbus, Germany (e-mail: lorenz-meyer@b-tu.de).}
		\thanks{J. G. Rueda-Escobedo is with the Basic Sciences Division, Engineering Faculty, National Autonomous University of Mexico, 04510 Mexico City, Mexico (e-mail:juan.rueda@fi.unam.edu).}
		\thanks{J. A. Moreno is with Electrical and Computer Engineering, Institute of Engineering, National Autonomous University of Mexico, 04510 Mexico City, Mexico, (e-mail:JMorenoP@ii.unam.mx).}
		\thanks{J. Schiffer is with Brandenburg University of Technology Cottbus-Senftenberg, 03046 Cottbus, Germany and Fraunhofer Research Institution for Energy Infrastructures and Geothermal Systems (IEG), 03046 Cottbus, Germany (e-mail: schiffer@b-tu.de).}
	}

	\maketitle
	
	\begin{abstract}
	While distributed parameter estimation has been extensively studied in the literature, little has been achieved in terms of robust analysis and tuning methods in the presence of disturbances. However, disturbances such as measurement noise and model mismatches occur in any real-world setting. Therefore, providing tuning methods with specific robustness guarantees would greatly benefit the practical application. To address these issues, we recast the error dynamics of a continuous-time version of the widely used consensus + innovations-based distributed parameter estimator to reflect the error dynamics induced by the classical gradient descent algorithm. This paves the way for the construction of a strong Lyapunov function. Based on this result, we derive linear matrix inequality-based tools for tuning the algorithm gains such that a guaranteed upper bound on the $\mathcal{L}_2$-gain with respect to parameter variations, measurement noise, and disturbances in the communication channels is achieved. An application example illustrates the efficiency of the method.
	\end{abstract}
	
	
	\IEEEpeerreviewmaketitle
	
	\section{Introduction}
	\subsection{Motivation and literature review}
	Distributed parameter estimation plays a major role in many modern engineering and scientific applications due to the large and increasing amount of sensors deployed in the field. Centralized processing, e.g., via a fusion center, has the drawback of introducing a single point of failure and potential congestion at certain communication links, which channel data to the fusion center, or violation of privacy considerations. This is especially the case for large and geographically dispersed networks \cite{kar_consensus_2013, xie_necessary_2018}. Consequently, distributed architectures for collaborative parameter estimation without a central entity have received an increasing amount of attention in the literature (see, e.g., \cite{kar_consensus_2013, xie_necessary_2018, papusha_collaborative_2014, chen_distributed_2014}). One promising approach for solving this problem is widely called \ac{c+i} \cite{kar_consensus_2013}. It is based on combining distributed consensus to reach global parameter agreement with local innovation utilizing measurements to update the local parameter estimation. Various extensions to this method have been proposed, many coming from a stochastic systems background \cite{zhang_distributed_2012, wang_distributed_2020, wang_decentralized_2021}. In these works, stationary system signals, i.e., constant regressors disturbed by noise and typically denoted as random measurement matrices, are assumed. This hinders the application to many cyber-physical and feedback systems, in which the system signals are not constant~\cite{xie_analysis_2018}. In a stream of papers by Xie et al., this assumption is dropped, and the algorithm is extended to the case of time-varying regressors \cite{xie_analysis_2018,xie_necessary_2018, xie_analysis_2020}. In \cite{xie_necessary_2018}, tracking error bounds are established, a performance analysis for the tracking error covariance matrix is performed, and recommendations for gain selection are given for the case of scalar outputs. 
	
	Starting with the papers \cite{papusha_collaborative_2014,chen_distributed_2014}, \ac{c+i}-type algorithms have been studied in a deterministic, continuous-time framework, naturally incorporating time-varying regressors.
	In those works, the notion of \ac{cpe}, an extension of the classical \ac{pe} condition \cite{anderson_exponential_1977}, is introduced. Recently, efforts have been made to extend the method to directed graphs. Results for strongly and weakly connected directed graphs are proposed in \cite{hu_adaptive_2021} and \cite{javed_excitation_2022}, respectively. 	
	In \cite{javed_excitation_2022}, the notion of \ac{uco} and the observability conservation under output feedback \cite{sastry_adaptive_1989} is utilized to establish the estimator's convergence. 
	
	Although extensive research has been conducted on \ac{c+i}-based distributed parameter estimation and several interesting extensions have been proposed in the literature, a notable gap remains regarding robust analysis and tuning conditions to robustify the algorithm in the presence of disturbances. Yet, in any real-world parameter estimation problem, disturbances in the form of measurement noise, perturbations in the communication channels, and model mismatches appear. As a result, tuning techniques with specific robustness guarantees would greatly benefit the practical application.   
	
	A key obstacle is that --- to the best of the author's knowledge --- no strong \ac{lf} is available thus far for a C+I-type distributed parameter estimator. Strong \acp{lf} provide means for calculating the convergence velocity and conducting robust analysis, including the study of input/output properties, such as \ac{iss} \cite{rueda-escobedo_strong_2021}. Moreover, a strong \ac{lf} can facilitate the derivation of tuning criteria for the algorithm's gains in the presence of disturbances~\cite{rueda-escobedo_l2-gain_2022}.  
	Recently, methods of constructing strong \acp{lf} have been proposed for the classical gradient descent algorithm in \cite{rueda-escobedo_strong_2021} and the consensus problem with single and double integrator dynamics in \cite{chowdhury_estimation_2018}. 
	\subsection{Contributions}
	Motivated by the aforementioned challenges and developments, we derive a strong \ac{lf} for a continuous-time \ac{c+i}-based distributed parameter estimator. Based on this result, we provide an $\mathcal{L}_2$-gain tuning method for the algorithm's gains in the presence of disturbances. More precisely, our main contributions are four-fold:
	\begin{itemize}
		\item We recast the error dynamics of a \ac{c+i}-based distributed parameter estimator in the form as induced by the classical gradient descent algorithm, paving the way for the subsequent robustness analysis.    
		\item By employing a scaled version of the observability Gramian as a strictifying term, we present a strong \ac{lf} for the \ac{c+i}-type distributed parameter estimator. To exemplify the application of this strong \ac{lf}, we derive a convergence bound for the error dynamics and the \ac{iss}-gain w.r.t. an affine disturbance acting on the error dynamics.
		\item By utilizing the derived strong \ac{lf}, we provide \ac{lmi} conditions for the upper and lower bounds on the algorithm's gradient descent gain. Those guarantee a maximum $\mathcal{L}_2$-gain of the global error dynamics w.r.t. a specified performance output. Building on this result, we provide an optimization routine to obtain the algorithm gains by minimizing the upper bound on the $\mathcal{L}_2$-gain, which we recall as the counterpart to the $H_\infty$-norm for \ac{ltv} and nonlinear systems.
		\item Finally, we illustrate the tuning method with an application example. 	
	\end{itemize}
	The remainder of the paper is organized as follows. In Section~\ref{sec:c+i}, a \ac{c+i}-based distributed parameter estimator is introduced and recast as a classical gradient descent algorithm for a nominal and a disturbed case. The main results of this paper are presented in Section~\ref{sec:main_result}. An application example is given in Section~\ref{sec:acc_ex}. Final remarks are given in Section~\ref{sec:conclusions}. The proofs of our results are provided in Section~\ref{sec:proofs}.
	\subsection*{Notation}
	The set of real numbers is denoted by $\mathbb{R}$, the $n$-dimensional Euclidean space by $\mathbb{R}^{n}$, and the set of real $n\times m$ matrices by $\mathbb{R}^{n\times m}$. Furthermore, we define the set $\mathbb{R}_{>0} \coloneqq \{x\in \mathbb{R} | x> 0\}  $. The symbol $\mrm{e}$ represents the Euler number. The vectors of ones and zeros with dimension $n$ are denoted by $\bm{1}_n$ and $\bm{0}_n$, respectively. The $n\times n$ identity matrix is indicated by $I_{n}$. Whenever the dimensions are clear from the context, by $0$, we mean a matrix of zeros. Otherwise, we denote a $n\times m$ matrix of zeros by $0_{n\times m}$. The Kronecker product is represented by $\otimes$, $\diag(A_1,\cdots,A_k)$ is used in a non-standard way and means a $k\,n \times k\,m$ matrix with diagonal entries $A_i\in \mathbb{R}^{n\times m}$, $i=\{1,\cdots,k\}$, and $\col(a_1,\cdots,a_k)$ denotes a vector obtained by stacking the vectors $a_i \in \mathbb{R}^{q}$, $i=\{1,\cdots,k\}$. For $A \in \mathbb{R}^{n \times n}$ symmetric, $A > 0 \ (A \ge 0)$ means that $A$ is symmetric positive (semi-)definite, with its largest and smallest eigenvalues denoted by $\lambda_\mrm{max}(A)$ and $\lambda_\mrm{min}(A)$, respectively. For $v \in  \mathbb{R}^{n}$, $\| v\|$ means $(v^\top v)^\frac{1}{2}$, and for $B \in  \mathbb{R}^{n \times m}$, $\|B\|$ denotes the induced matrix norm of $B$, defined as $\mrm{sup}_{\|v\|=1}\|Bv\|$. The space of vector-valued functions with $n$ components that are square-integrable are represented by $\mathcal{L}_2(n,[0,\infty))$. When there is no risk of misunderstanding, we use the notation $\mathcal{L}_2$ to represent this space.
	The $\mathcal{L}_2$-norm of $u(t)$ is defined as $\|u(t)\|^2_{\mathcal{L}_2} = \underset{t \to \infty}{\mrm{lim}} \int_{0}^{t} u^\top(s)u(s) \mrm{d}s$. A system $H : \mathcal{L}_2(n, [0,\infty)) \to\mathcal{L}_2(m, [0,\infty))$ is considered to have a finite $\mathcal{L}_2$-gain if there exist non-negative constants $\gamma$ and $b$, such that \cite{van_der_schaft_l2-gain_2017} 
	\begin{equation}
		\label{eq:L2_gain_system}
		\|H(u(t))\|_{\mathcal{L}_2} \le \gamma \|u(t)\|_{\mathcal{L}_2} + b, \ \forall \ u(t) \in \mathcal{L}_2(n,[0,\infty)) .
	\end{equation} 
	Moreover, the $\mathcal{L}_2$-gain of $H$ is defined as
	$$\gamma^* \coloneqq \mrm{inf} \ \{\gamma  |  \exists b  \text{ such that \eqref{eq:L2_gain_system} holds} \}.$$
	\section{C+I-based distributed parameter estimation}
	\label{sec:c+i}
	In this paper, we consider a sensor network with $n$ sensors, where each one is taken as a node in a multi-agent system. We represent the communication structure among the agents by an undirected, time-dependent graph $\mathcal{G}(t) = \left(\mathcal{V}, \mathcal{E}(t)\right)$, where $\mathcal{V}= \{1,2,...,n\}$ is the set of nodes, and $\mathcal{E}(t)\subseteq\mathcal{V} \times \mathcal{V}$ is the set of edges at time $t$, and the cardinality of $\mathcal{E}(t)$ is denoted by $n_e(t)$. Thus, we allow the relationship between neighbors to change over time, e.g., to model communication failures. By assigning an arbitrary orientation to the graph  $\mathcal{G}(t)$ at time $t$, the oriented incidence matrix $D(t) \in \mathbb{R}^{n \times n_\mrm{e}(t) }$ is defined element-wise as $D_{ij}(t) = 1$, if $i$ is the source of the $l$th egde, $D_{ij}(t) = 1$, if $i$ is the sink of the $l$th egde, and $D_{ij}(t) = 0$ otherwise.	
	The Laplacian matrix $L(t) \in  \mathbb{R}^{n\times n}$ of the graph $\mathcal{G}(t)$ at time $t$ is defined as $L(t) = D(t)D^\top(t)$.
	
	\subsection{Nominal system}
	At each agent, we consider the following measurement model in the form of a \ac{lre}~\cite{narendra_stable_2012}
	\begin{equation}
		\label{eq:local_LRE}			
		y_i(t) = C_i(t) \theta,
	\end{equation}
	where $i \in \mathcal{V}$, $y_i \in \mathbb{R}^{N_\mrm{y}}$ denotes the local output of the $i$th agent, $C_i(t)\in \mathbb{R}^{N_\mrm{y}\times N}$ the local regressor of the $i$th agent, and $\theta \in \mathbb{R}^N$ the global parameter vector.
	To estimate consistent parameters $\theta$ at all agents and building on the algorithms proposed in, e.g., \cite{papusha_collaborative_2014, chen_distributed_2014}, we utilize the following \ac{c+i}-type algorithm for each agent to update the local estimate $\hat{\theta}_i(t) \in \mathbb{R}^N$ of $\theta$ at agent $i$
	\begin{equation}
		\begin{aligned}
			\label{eq:C+I_alg}
			\dot{\hat{\theta}}_i(t) = &-\alpha\Gamma_i \sum_{j\in\mathcal{N}_i(t)}(\hat{\theta}_i(t)-\hat{\theta}_j(t)) \\ &- \Gamma_i C_i^\top(t)(C_i(t) \hat{\theta}_i(t)-y_i(t)),
		\end{aligned}	
	\end{equation}
	where the subscript $i$ denotes the variables associated to the $i$th agent, $\mathcal{N}_i(t)$ denotes the set of neighbors of the agent at time $t$, $\Gamma_i \in \mathbb{R}^{N\times N}$ is the gradient descent gain matrix, $\alpha \in \mathbb{R}$ is the additional consensus gain,.
	\begin{remark}
		By choosing $\alpha =1$ and $\Gamma_i = \gamma_\mrm{p}I_{N}$ with $\gamma_\mrm{p}>0$ in \eqref{eq:C+I_alg}, we recover the algorithm of \cite{papusha_collaborative_2014, chen_distributed_2014}.
	\end{remark}
	To write the \ac{c+i}-type algorithm in a compact form, we stack each agent's estimate as 
	\begin{align*}
		\hat{\Th}(t) = \col\left(\hat{\theta}_1^\top(t), \cdots, \hat{\theta}_n^\top(t)\right),	
	\end{align*}
	and write the dynamics of $\hat{\Th}$ as 
	\begin{align}
		\dot{\hat{\Th}}(t) =&-\alpha \bar{\Gamma}\,\bar{L}(t)\hat{\Th}(t)-\bar{\Gamma}\,\bar{C}^\top(t)\big(\bar{C}(t)\hat{\Th}(t)-\bar{y}_1(t)\big), \label{eq:dot_hat_theta}
	\end{align}
	where $\bar{\Gamma} = \diag(\Gamma_1,\cdots, \Gamma_n)$, $\bar{L}(t) = L(t) \otimes I_{N},$ and
	\begin{align*}
		\bar{C}(t) &= \diag(C_1^\top(t),\cdots, C_n^\top(t)),\\
		\bar{y}_1(t) &= \col(y_1(t), ..., y_n(t)).
	\end{align*}
	We recall that $L(t) \in  \mathbb{R}^{n\times n}$ is the Laplacian matrix of the graph $\mathcal{G}(t)$ at time $t$. 
	
	We make the following assumption on the algorithm gains, which can always be fulfilled by proper design.
	\begin{assumption}
		\label{ass:2}
		For each agent, its gradient descent gain matrix $\Gamma_i$ in \eqref{eq:C+I_alg} is symmetric and positive definite, so that $\bar{\Gamma}$ in \eqref{eq:dot_hat_theta} is symmetric and positive definite, with norm $r_1:=\|\bar{\Gamma}\|$. Moreover, the additional consensus gain $\alpha$ in \eqref{eq:C+I_alg} is chosen to be strictly positive.
	\end{assumption}	
	Additionally, we define the estimated output generated by the system \eqref{eq:dot_hat_theta} as
	\begin{equation}
		\begin{aligned}
			\yhat(t) &= \begin{bmatrix} \hat{y}_a(t) \\ \hat{y}_b(t)\end{bmatrix} = \bar{\Lambda}(t)\hat{\Th}(t),\\
			\bar{\Lambda}(t) &= \begin{bmatrix}\bar{C}(t) \\ \sqrt{\alpha} D^\top(t) \otimes I_{N}\end{bmatrix}\in \mathbb{R}^{\left(nN_\mrm{y}+Nn_{\mrm{e}(t)}\right)\times nN},
		\end{aligned}\label{eq:y_hat}
	\end{equation}
	where $D(t) \in \mathbb{R}^{n \times n_\mrm{e}(t) }$ denotes the oriented incidence matrix of the graph $\mathcal{G}(t)$ at time $t$. Here, ${\hat{y}_a(t)} \in  \mathbb{R}^{N_\mrm{y}n}$ corresponds to the stacked estimated measurements of each agent, and ${\hat{y}_b(t)}\in  \mathbb{R}^{N n_\mrm{e}(t)}$ to the difference of the estimations across each edge scaled by $\sqrt{\alpha}$. By defining the estimated output in this way, we preserve the linear regression structure introduced in \eqref{eq:local_LRE}, i.e., $\yhat(t) = \la(t) \hat{\Th}(t)$, for the network setting. This allows us to recast the error dynamics of a \ac{c+i}-based distributed parameter estimator to reflect the error dynamics induced by the classical gradient descent algorithm in the following. 
	
	The true parameter vector of the global system is
	\begin{equation*}
		{\Th} =  \bm{1_n}\otimes \theta .
	\end{equation*}
	Thus, to recast the error dynamics, we introduce the output 
	\begin{equation}
		\label{eq:output}
		\y(t) =    {\la(t)} {\Th}(t) = \begin{bmatrix}
			{\bar{y}_a(t)} \\ \bm{0}_{N n_\mrm{e}(t)} \end{bmatrix}, 
	\end{equation}
	where ${\bar{y}_b(t)} = \bm{0}_{N n_\mrm{e}(t)}$ follows from the relation
	\begin{align*}
		(D^\top(t) \otimes I_{N}) (\bm{1_n}\otimes\theta) =  \underbrace{(D^\top(t)\bm{1}_n)}_{=\bm{0}_{n_\mrm{e}(t)}} \otimes(I_{N}\theta).
	\end{align*}
	To derive the relation above, we used the mixed product property of the Kronecker product and the fact that there is one element equal to $+1$, one element equal to $-1$, and all other zeros in each row of $D^\top(t)$ \cite{bullo_lectures_2020}. Hence, we can rewrite \eqref{eq:dot_hat_theta} in the form of a classical gradient descent algorithm, i.e.,
	\begin{equation}
		\label{eq:dot_hat_theta2}
		\dot{\hat{\Th}}(t) = -\bar{\Gamma}\la^\top(t) \left(\la(t) {\hat {\Th}} (t) - {\y(t)}\right)	.
	\end{equation}
	\begin{remark}
		As \eqref{eq:dot_hat_theta2} resembles the form of a gradient descent algorithm, it is immanent to assert that it minimizes the instantaneous quadratic cost function $$J({\hat{\Th}}(t)) =   \frac{1}{2} \left(\y(t) -\yhat(t) \right)^\top\left(\y(t) -\yhat(t)  \right) .$$
		Hence, \eqref{eq:dot_hat_theta2} can be recovered from the gradient flow (cf. \cite{papusha_collaborative_2014, javed_excitation_2022}), i.e.,
		\begin{align*}
			\dot{\hat{\Th}}(t) = -\bar{\Gamma} \frac{\partial}{\partial {\hat{\Th}}(t)}J({\hat{\Th}}(t)).
		\end{align*}
	\end{remark}
	By defining the estimation error as $\tilde{\Th}(t)  \coloneqq \hat{\Th}(t)-{\Th},$ and the output error as $\ytil(t) \coloneqq \yhat(t)-\y(t)  ,$ the error system corresponding to \eqref{eq:dot_hat_theta2} can be rewritten as
	\begin{equation}
		\begin{split}
			\label{eq:err_dyn2}
			\dot{\tilde{\Th}}(t) &=  -\bar{\Gamma}\la^\top(t) \la(t) {\tilde{\Th}} (t)	, \\
			\ytil(t) &=  \la(t) \tilde{\Th}(t) .
		\end{split}
	\end{equation}
	This error system is the foundation for deriving the strong \ac{lf} in Section~\ref{sec:sLF} and, hence, proving \ac{guas} (or equivalently Uniform Exponential Stability) of the origin. It also facilitates further analysis, such as calculating the convergence velocity and the \ac{iss}-gain with respect to an affine disturbance.

	\section{Main results}
	\label{sec:main_result}
	To present our main results, we make the following assumptions on the \ac{lre} \eqref{eq:local_LRE} and the graph $\mathcal{G}(t)$ representing the communication structure among agents.
	\begin{assumption} \
		\label{ass:1}
		\begin{enumerate}[leftmargin=1.1cm, label=({A2.}{{\arabic*}})]
			\item The \ac{cpe} condition is fulfilled, i.e., there exist positive constants $T$, $\cPEu \ge\cPEl > 0$, all independent of $t$, such that for all $t\ge T$ it holds that
			\begin{equation}
				\label{eq:cPE}
				\cPEu I_{N} \ge \int_{t-T}^{t}\sum_{i=1}^{n}{C_i}^\top(s){C_i}(s)\mathrm{d}s \ge \cPEl I_{N}  .
			\end{equation}
			\item All the local regressors are uniformly upper-bounded in the norm, with upper bound given by $r_2>0$, i.e.,
			\begin{align*}
				\|C_i(t)^\top C_i(t)\| \le  r_2 \quad \forall\,t\geq 0 \quad \forall i=1,\cdots,n.
			\end{align*}
			\item The graph $\mathcal{G}(t)$ is connected on average, meaning that $\int_{t-T}^{t}L(s)\mathrm{d}s$, with $T$ as in A2.1, has only one zero eigenvalue $\lambda_1$, and all remaining ones are positive and accept a lower bound denoted by the constant $0<\underline{\lambda} \le \lambda_i(t) , \  \forall i \ge 2$. Furthermore, the time-varying Laplacian $L(t)$ in \eqref{eq:dot_hat_theta} is uniformly upper-bounded in the norm, with upper bound $r_3>0$, i.e.,
			\begin{align*}
				\|L(t)\| \le r_3 \quad \forall\,t\geq 0.	
			\end{align*}
		\end{enumerate} 
	\end{assumption}
	Assumption~\ref{ass:1} follows standard practice and is commonly imposed in the literature (see, e.g., \cite{anderson_exponential_1977,rueda-escobedo_strong_2021} for A2.2 and \cite{papusha_collaborative_2014, chen_distributed_2014} for A2.1 and A2.3). Furthermore, the assumption on the required excitation A2.1 is the same for a distributed and a centralized algorithm. This is recognized by stacking the local \acp{lre} \eqref{eq:local_LRE}, to obtain a global \ac{lre} in the form
	\begin{equation*}
		\y^\mrm{c}(t) = \col\left(C_1(t), \cdots,\, C_n(t)\right)\theta = C^\mrm{c}(t)\theta.
	\end{equation*}
	The parameters can be estimated, e.g., by a centralized gradient descent algorithm. In such case, the error dynamics is given by
	\begin{equation}
		\label{eq:centr_gradient_descent}
		\frac{\mrm{d}}{\mrm{d}t} (\hat{\theta}(t)-\theta)= - \Gamma^\mrm{c}(C^\mrm{c})^\top(t)C^\mrm{c}(t)\left(\hat{\theta}(t)-\theta\right),
	\end{equation}
	where $\Gamma^\mrm{c}\in \mathbb{R}^{N\times N}$ denotes the centralized gradient descent gain.
	A necessary and sufficient condition for \ac{guas} of the origin of \eqref{eq:centr_gradient_descent} is \ac{pe} of the global regressor $C^\mrm{c}$ \cite{anderson_exponential_1977}, i.e.,
	\begin{equation*}
		\cPEu^\mrm{c} I_{N} \ge \int_{t-T^\mrm{c}}^{t}(C^\mrm{c})^\top C^\mrm{c}\mathrm{d}s \ge \cPEl^\mrm{c} I_{N} \qquad \forall t\ge T^\mrm{c},
	\end{equation*}
	where $\cPEu^\mrm{c} \ge\cPEl^\mrm{c} > 0$ and $T^\mrm{c} > 0$ are constants independent of $t$. 
	This condition is evidently equivalent to \ac{cpe} \eqref{eq:cPE}. 
	
	For ease of presentation, we present the proofs of all the following results in Section~\ref{sec:proofs}. Furthermore, we define $r_4 \coloneqq r_2+\alpha r_3$, with $r_2$ and $r_3$ defined in Assumption \ref{ass:1}.
	
	\subsection{Strong Lyapunov function for C+I-based distributed parameter estimation}
	\label{sec:sLF}
	To streamline our results, we define the system  
	\begin{equation}
		\label{eq:err_dyn2_OI}
		\begin{split}
			\dot{\tilde{\Th}}^\mrm{OI}(t) &=   \bm{0}_{nN}	, 	
			\\
			\ytil^\mrm{OI}(t) &=  \la(t) \tilde{\Th}^\mrm{OI}(t) ,
		\end{split}
	\end{equation}
	with output as in \eqref{eq:output} and obtained from the error system \eqref{eq:err_dyn2} by a linear output injection. 
	To construct a strong \ac{lf} for the error system \eqref{eq:err_dyn2}, we first establish the \ac{uco} of the system \eqref{eq:err_dyn2_OI}. Using the observability conservation by output injection \cite{zhang_observability_2015}, we would later prove the \ac{uco} of the error system \eqref{eq:err_dyn2}. Finally, this allows us to present a strong \ac{lf}. 
	
	To proceed, we denote the observability Gramian of the system \eqref{eq:err_dyn2_OI} as
	\begin{equation}
		\label{eq:gram_OI}
		M^\mrm{OI}(t,t-T) = \int_{t-T}^{t}\la^\top(s)\la(s)\mathrm{d}s .
	\end{equation}
	We have the following result about the \ac{uco} of the system \eqref{eq:err_dyn2_OI}.
	\begin{lemma}
		\label{lemma:obs_gram_OI}
		Under Assumptions~\ref{ass:2}~and~\ref{ass:1}, \ac{cpe} \eqref{eq:cPE} implies the \ac{uco} of the system \eqref{eq:err_dyn2_OI}, with the bounds on $M^\mrm{OI}$ in \eqref{eq:gram_OI} as follows:
		\begin{equation*}
			\MOIu  I_{nN}\ge	M^\mrm{OI}(t,t-T) \ge \MOIl I_{nN},
		\end{equation*}
		with
		\begin{equation*}
			\begin{split}
				\MOIu &= Tr_4   , 
				\\ \MOIl &= \underset{\forall \left\| a(t)\right\|^2 \le 1}{\mrm{min}} \mrm{max} \bigg\{ \alpha  \underline{\lambda} ( 1-\|a(t)\|^2 ), \\ &\frac{\cPEl}{n} \|a(t)\|^2 -    2 r_2 T  \sqrt{ \left\| a(t)\right\|^2 (1-\left\|a(t)\right \|^2)} \bigg\} >0 ,
		\end{split}\end{equation*}
		and where $a(t) \coloneqq [a_1(t), ..., a_N(t)]^\top \in \mathbb{R}^N $ is a vector with components $a_i(t) \ge 0, \ \forall t \in \mathbb{R}, \ \forall i \in  \mathcal{V}$ defined in \eqref{eq:linear_combination_EV_L}.
	\end{lemma}
	Moreover, the observability Gramian of the system \eqref{eq:err_dyn2} corresponds to
	\begin{equation}
		\label{eq:gram}
		M(t,t-T) = \int_{t-T}^{t}\bar{\Phi}^\top(s,t)\la^\top(s)\la(s)\bar{\Phi}(s,t)\mathrm{d}s ,
	\end{equation}
	where $\bar{\Phi}(s,t)$ is the state transition matrix of the system \eqref{eq:err_dyn2}.
	Adapting \cite[Lemma 2]{zhang_observability_2015} to the systems \eqref{eq:err_dyn2_OI} and \eqref{eq:err_dyn2}, we present the following result. 
	\begin{lemma}
		\label{lemma:obs_gram}
		Under Assumptions~\ref{ass:2}~and~\ref{ass:1}, the system \eqref{eq:err_dyn2} is \ac{uco} and its observability Gramian \eqref{eq:gram} accepts the following bounds:
		\begin{equation*}
			\Mu I_{nN}  \ge	M(t,t-T) \ge \Ml I_{nN},
		\end{equation*}
		with
		\begin{equation}
			\label{eq:iota_3}
			\begin{split}
				\Mu &= \left(\sqrt{\MOIu-\MOIl+\varphi_2\MOIl}+\sqrt{\MOIu}\right)^2, 
				\\ \Ml &= \begin{cases}
					\left(\frac{\sqrt{\MOIu-\MOIl+\varphi_1\MOIl}-\sqrt{\MOIu}}{\varphi_1-1}\right)^2 >0& \text{if $\varphi_1 \ne 1$} ,\\
					\frac{\MOIl}{4\MOIu}>0& \text{if $\varphi_1 = 1,$} 
				\end{cases}
			\end{split}
		\end{equation}
		where
		\begin{gather}
			\varphi_1 = \frac{1}{2} r_1^2 \MOIu^2 \quad\text{and}\quad \varphi_2 = \frac{1}{4} \left(\mrm{e}^{2r_1\MOIu} -1 \right)-\frac{1}{2}r_1\MOIu. \label{eq:varphi}
		\end{gather}
	\end{lemma}
	\begin{remark}
		To obtain sharper bounds, the results of \cite[Lemma 4]{zhang_observability_2015} can be used with $\varphi_1$ and $\varphi_2$ as given in \eqref{eq:varphi}. 
	\end{remark}
	Knowing that \eqref{eq:err_dyn2} is \ac{uco}, we use a scaled version of its observability Gramian as a strictifying term in the strong \ac{lf} presented in the following theorem.
	\begin{theorem}
		\label{th:stong_lya}
		Consider the error system \eqref{eq:err_dyn2} of the \ac{c+i}-based distributed parameter estimator \eqref{eq:C+I_alg} together with Assumptions~\ref{ass:2}~and~\ref{ass:1}. 
		Then, the origin of \eqref{eq:err_dyn2} is globally uniformly asymptotically stable, and the quadratic function 
		\begin{equation}
			V(\tilde{\Th}, t) = \tilde{\Th}^\top(t)P(t)\tilde{\Th}(t),
			\label{eq:lya_fun}
		\end{equation}
		with 
		\begin{equation*}
			P(t) = \frac{T}{2}\bar{\Gamma}^{-1} + \int_{t-T}^{t}(s-t+T)\bar{\Phi}^\top(s,t)\la^\top(s)\la(s)\bar{\Phi}(s,t)\mathrm{d}s
		\end{equation*}
		is a strong \ac{lf} for the system \eqref{eq:err_dyn2}.
	\end{theorem}
		\begin{remark}
			Extending the result of \cite{rueda-escobedo_strong_2021} for the classical gradient descent parameter estimator for the system \eqref{eq:err_dyn2}, it is possible to use a scaled version of the observability Gramian of the system \eqref{eq:err_dyn2_OI} as a strictifying term for the construction of a strong \ac{lf}. 
%
			The reason to use \eqref{eq:lya_fun} instead in this work is that stricter bounds can be obtained.
		\end{remark}

	This strong \ac{lf} \eqref{eq:lya_fun} can be used to compute bounds on the rate of convergence of the error dynamics of the system \eqref{eq:err_dyn2} as is shown in the following corollary.
	\begin{corollary}
		\label{corollary:convergence_bounds}
		Consider the error system \eqref{eq:err_dyn2} of the \ac{c+i}-based distributed parameter estimator \eqref{eq:C+I_alg} together with Assumptions~\ref{ass:2}~and~\ref{ass:1}. Under these conditions, the trajectories of the error system \eqref{eq:err_dyn2} are bounded by
		\begin{equation}
			\label{eq:trajecty_bound}
			\|\tilde{\Th}(t)\| \le  	\sqrt{\frac{\kappa_1}{\kappa_2}} \|\tilde{\Th}(t_0)\|\,\mrm{e}^{-\frac{ \Ml}{2\kappa_1}(t-t_0)},
		\end{equation}
		with 
		\begin{equation}
			\label{eq:kappa}
			\kappa_1 =	\left( \frac{T}{2}\lambda_\mrm{max}(\bar{\Gamma}^{-1}) + T\Mu\right), \quad \kappa_2 = \frac{T}{2} \lambda_\mrm{min}(\bar{\Gamma}^{-1})  .
		\end{equation}
	\end{corollary}
	The applicability of the strong \ac{lf} in \eqref{eq:lya_fun} does not end here since it can be used to investigate the \ac{iss}-gains of the system \eqref{eq:err_dyn2} (see \cite[Chapter 4.9]{khalil_nonlinear_2002}). For this, we consider an affine disturbance $\delta_\mrm{ISS}(t)$ acting on the error dynamics as 
	\begin{equation}
		\label{eq:disturbed_system_ISS}
		\dot{\tilde{\Th}}^\mrm{ISS}(t) = -\bar{\Gamma}\la^\top(t)\la(t){\tilde{\Th}}^{\mrm{ISS}}(t) + \delta_\mrm{ISS}(t)	.
	\end{equation}
	\begin{corollary}
		\label{corollary:ISS_bounds}
		Consider the error dynamics \eqref{eq:disturbed_system_ISS} of the \ac{c+i}-based distributed parameter estimator \eqref{eq:C+I_alg} under the affine disturbance $\delta_\mrm{ISS}(t)$ and assume that Assumptions~\ref{ass:2}~and~\ref{ass:1} hold. Then, an upper bound for the \ac{iss}-gain from the disturbance $\delta_\mrm{ISS}(t)$ to the state ${\tilde{\Th}}^\mrm{ISS}(t)$ is
		\begin{equation} \label{eq:ISS_gain}
			\gamma_{\mrm{ISS}} \le  \frac{2\kappa_1}{\Ml}\sqrt{\frac{\kappa_1}{\kappa_2}},
		\end{equation} 
		with $\kappa_1$ and $\kappa_2$ given in \eqref{eq:kappa}.
	\end{corollary}
	As seen from Corollaries~\ref{corollary:convergence_bounds} and \ref{corollary:ISS_bounds}, the chosen gains $\alpha$ and $\bar{\Gamma}$ directly influence the convergence velocity, the ISS-gain and, thus, the performance of the algorithm. This observation motivates the development of a robust tuning method in the next section. 
	
	\subsection{$\mathcal{L}_2$-gain tuning in the presence of disturbances}
	\label{sec:tuning}
	The \ac{guas} of the origin of \eqref{eq:err_dyn2} is independent of the choice of the gains $\bar{\Gamma}>0$ and $\alpha>0$, as shown in Theorem~\ref{th:stong_lya}, and only depends on Assumptions \ref{ass:2} and \ref{ass:1} However, the performance of the \ac{c+i}-based distributed parameter estimator \eqref{eq:C+I_alg} is strongly affected by those gains. For different conditions, depending on the signals excitation level, measurement noise, communication disturbances, and possible parameter variations, an intelligent choice of gains should improve the performance. Under these conditions, a trade-off between convergence speed and disturbance attenuation is desirable.
	
	To derive tuning conditions for the gains $\bar{\Gamma}$ and $\alpha$, which guarantee an upper bound for the $\mathcal{L}_2$-gain of the global error dynamics w.r.t. a specified performance output, we introduce a (deterministic) disturbance $\delta(t)\in \mathbb{R}^{r}$ 
	that models parameter variation, measurement noise, and disturbances in the communication among agents. This disturbance is projected in the parameter dynamics with the aid of the matrix $\D_1\in  \mathbb{R}^{N\times r}$ and in the following manner:
	\begin{equation*}
		\dot{\theta}^\mrm{d}(t) = \D_1 \delta(t). 
	\end{equation*}
	The true parameter vector of the global system, hence, shows the following dynamics 
	\begin{equation*}
		\dot{\Th}^\mrm{d}(t) =  \bm{1}_n\otimes \dot{\theta}^\mrm{d}(t) = \Db1\,\delta(t), \quad \Db1=\bm{1}_n\otimes \D_1.
	\end{equation*}
	In the case of the output, $\delta(t)$ is projected as follows:
	\begin{equation*}
		\y^\mrm{d}(t)= \la(t) \Th^\mrm{d}(t) +  \Db2 \delta(t),
	\end{equation*}
	with $\Db2 \in  \mathbb{R}^{(N_\mrm{y}n+Nn_\mrm{e}(t))\times r}$. Recall that the output defined in \eqref{eq:output} stacks the measurements of each agent and a scaled version of the difference of the estimations across each edge (see \eqref{eq:y_hat}). Thus, $\Db2$ may project different parts of $\delta(t)$ to these two output components. By defining the disturbed global estimation error as $\tilde{\Th}^\mrm{d}(t) \coloneqq	 \hat{\Th}^\mrm{d}(t)  - {\Th}^\mrm{d}(t)$, the disturbed estimation error dynamics results in
	\begin{equation}
		\label{eq:err_dyn_dis}
		\begin{split}
			\dot{\tilde{\Th}}^\mrm{d}(t) &=	\dot{\hat{\Th}}^\mrm{d}(t)-\dot{{\Th}}^\mrm{d}(t) \\&= -\bar{\Gamma}\la^\top(t) \la(t){\tilde{\Th}}^\mrm{d} (t)	+ \left(\bar{\Gamma}\la^\top(t)\Db2-\Db1\right){\delta} (t),
		\end{split}
	\end{equation}
	where $\dot{\hat{\Th}}^\mrm{d}(t)$ is obtained from \eqref{eq:dot_hat_theta2} using the disturbed output $\y^\mrm{d}(t)$ instead of the nominal output $\y(t)$. The system \eqref{eq:err_dyn_dis} is a perturbed version of the nominal estimation error dynamics \eqref{eq:err_dyn2}. To investigate the robustness of \eqref{eq:err_dyn2} in terms of its $\mathcal{L}_2$-gain w.r.t. the perturbations, i.e., parameter variations, measurement noise, and disturbances in the communication among agents, we define the performance output
	\begin{equation}
		\label{eq:z}
		z(t) = Q\tilde{\Th}^\mrm{d}(t) + W \delta(t),
	\end{equation}
	with $z(t)\in \mathbb{R}^p $, $Q \in \mathbb{R}^{p\times nN}$ and $W \in \mathbb{R}^{p\times r}$. We assume the standard restrictions $W^\top W= I_{r}$ and $Q^\top W = 0$. Hence, we have 
	\begin{equation}
		\label{eq:zTz}
		z^\top(t) z(t) = (\tilde{\Th}^\mrm{d})^\top(t) Q^\top Q\tilde{\Th}^\mrm{d}(t) +  \delta(t)^\top\delta(t).
	\end{equation}

	Inspired by the methodology introduced in \cite{rueda-escobedo_l2-gain_2022}, we present the following result to quantify the effect of the disturbance $\delta(t)$ on the performance output $z(t)$ as defined in \eqref{eq:z} via the $\mathcal{L}_2$-norm.	
	For this, we formulate the following problem. 
	\begin{problem}
		\label{prob:1}
		Fix the constants $\gamma>0$, $\alpha>0$, and $c_2>0$ and find parameters $c_1 >0$, $\gamma_1>0$, and $\gamma_2 >0$, such that 
		\begin{equation}
			\label{eq:LMI_tuning_th1}
			\begin{split}
				\gamma_1 - \gamma_2 &\ge 0, \\
				\begin{bmatrix}
					\Phi_{11} & \Phi_{12}\\ \Phi_{12}^\top & -\gamma_2 I_{nN}
				\end{bmatrix} &\le0,
			\end{split}
		\end{equation}
		where $\Phi_{11}$ and $\Phi_{12}$ are defined in \eqref{eq:Phi}.
	\end{problem}
		The chosen parameters $\gamma$, $\gamma_1$, $\gamma_2$, $\alpha$, $c_1$, and $c_2$, as well as the excitation parameters $\Mu$, $\Ml$, $r_4$ and $T$, influence the feasibility of the \ac{lmi} \eqref{eq:LMI_tuning_th1}. However, by choosing large enough values for the parameters $\gamma$, $\gamma_2$, $c_1$, and $c_2$ it is possible to obtain a solution that renders the \ac{lmi} \eqref{eq:LMI_tuning_th1} feasible. With the known parameters $\gamma_1$ and $\gamma_2$, the gain $\bar{\Gamma}$ can be selected as given in the following theorem. 
	\begin{theorem}
		\label{th:tuning}
		Consider the system \eqref{eq:err_dyn_dis} together with Assumption~\ref{ass:1}, as well as input $\delta(t)$ and output $z(t)$ given in \eqref{eq:z}. Fix the constants $\gamma>0$ and $\alpha>0$, and choose $c_2>0$, such that 
		\begin{equation}
			\label{eq:c2_LMI}
			\begin{bmatrix}
				\frac{c_2\Ml}{2} & -Q^\top \\ -Q & I_p
			\end{bmatrix} >0.
		\end{equation}
		Suppose that Problem~\ref{prob:1} is feasible for $\gamma>0$, $\gamma_1>0$, $\gamma_2>0$, $\alpha>0$, $c_1>0$, and $c_2>0$. Choose $\bar{\Gamma}$, such that  
		\begin{equation}
			\label{eq:LMI_Gamma}
			\sqrt{\gamma_1} I_{nN} \ge \bar{\Gamma} \ge  	\sqrt{\gamma_2} I_{nN },
		\end{equation}
		then the system \eqref{eq:err_dyn_dis} with output \eqref{eq:z} has an $\mathcal{L}_2$-gain of at most $\sqrt{\gamma}$.
	\end{theorem}
	\begin{figure*}[!t]
		\begin{equation}
			\scriptsize
			\label{eq:Phi}
			\begin{aligned}
				\Phi_{11}&= 	\begin{bmatrix}\phi_{11}  & 0 & 0 & 0 \\ 
					0 & -\frac{c_2\Ml}{2}I_{nN} & Q^\top & 0 \\
					0 & Q & -I_{p} & W \\
					0 & 0 & W^\top  &  \phi_{44}  \end{bmatrix} \ , \quad 
				\Phi_{11}^\mrm{OE}= 	\begin{bmatrix}\phi_{11}  & 0 & (Q^\mrm{OE})^\top & 0 \\ 
					0 & -\frac{c_2\Ml}{2}I_{nN} &  0& 0 \\
					Q^\mrm{OE} & 0 & -I_{p} & W^\mrm{OE} \\
					0 & 0 & (W^\mrm{OE})^\top  &  \phi_{44}  \end{bmatrix} \ , \quad 
				\Phi_{12}= 	\begin{bmatrix}0 \\ 0 \\ 0 \\ -\frac{ 2c_1 }{\sqrt{c_2 \Ml}} \Db1^\top\end{bmatrix} \ , \\
				\phi_{11} &=    \left(-c_1 + c_2T \right)I_{\left(nN_\mrm{y}+Nn_{\mrm{e}(t)}\right)}   \ ,  \quad \phi_{44} = \frac{8c_2\Mu^2}{\Ml} \Db1^\top \Db1 +	 \left(c_1 +  \frac{8c_2\Mu^2r_4\gamma_1}{\Ml}\right) \Db2^\top \Db2   -\gamma I_{r}   .
			\end{aligned}
		\end{equation}
		\hrule 
	\end{figure*}
 By leaving $\gamma>0$ as a decision variable and combining the \ac{lmi} \eqref{eq:LMI_tuning_th1} with a minimization over $\gamma$, we can formulate a \ac{sdp}. From the solution of the \ac{sdp}, we can determine the gain $\bar{\Gamma}$ utilizing \eqref{eq:LMI_Gamma}. This solution depends on the chosen additional consensus gain $\alpha$. Thus, we combine the \ac{sdp} with a nonlinear optimization routine in the scalar decision variable $\alpha$. This idea gives rise to the following corollary.  
	\begin{corollary}
		\label{co:tuning}
		Consider the system \eqref{eq:err_dyn_dis} together with Assumption~\ref{ass:1} and input $\delta(t)$ and output $z(t)$ given in \eqref{eq:z}. Choose $\alpha>0$ as
		\begin{equation*}
			\underset{\alpha \in \mathbb{R}_{>0}}{\mrm{arg \ max}} \ \frac{\Ml(\alpha)}{\Mu(\alpha)},
		\end{equation*}
		with $\Ml$, and $\Mu$ defined in \eqref{eq:iota_3}. Choose $c_2>0$, such that \eqref{eq:c2_LMI} is satisfied.
		Solve the following \ac{sdp} in the scalar decision variables ($c_1, \gamma_1, \gamma_2, \gamma$):
		\begin{equation}
			\label{eq:sdp1}
			\begin{aligned}
				&\mrm{minimize} \ \ \gamma&\\
				&\mrm{subject \ to:}	 	& \gamma >0, \ \gamma_1 - \gamma_2 \ge 0, \\
				&				&\begin{bmatrix}
					\Phi_{11} & \Phi_{12}\\ \Phi_{12}^\top & -\gamma_2 I_{nN}
				\end{bmatrix} &\le0,
			\end{aligned}
		\end{equation}
		where $\Phi_{11}$ and $\Phi_{12}$ are defined in \eqref{eq:Phi}. If the \ac{sdp} \eqref{eq:sdp1} is feasible, choose $\bar{\Gamma}$, such that 
		\begin{equation*}
			\sqrt{\gamma_1} I_{nN } \ge \bar{\Gamma} \ge  	\sqrt{\gamma_2} I_{nN }.
		\end{equation*}
		In such case, the system \eqref{eq:err_dyn_dis} with output \eqref{eq:z} has an $\mathcal{L}_2$-gain of at most $\sqrt{\gamma}$.
	\end{corollary}
	
	Following \cite[Corollary 2]{rueda-escobedo_l2-gain_2022}, it is straightforward to extend the result of Corollary~\ref{co:tuning} to the modified output 
	\begin{equation}
		\label{eq:z_OE}
		z^\mrm{OE}(t) = Q^\mrm{OE}\la(t)\tilde{\Th}^\mrm{d}(t) + W^\mrm{OE} \delta(t),
	\end{equation} 
	with $Q^\mrm{OE} \in \mathbb{R}^{p\times \left(nN_\mrm{y}+Nn_{\mrm{e}(t)}\right)}$, $W^\mrm{OE} \in \mathbb{R}^{p\times r}$, and restrictions $(W^\mrm{OE})^\top W^\mrm{OE}= I_{r}$ and $(Q^\mrm{OE})^\top W^\mrm{OE} = 0$.
	Here, only the states reflected in the output error $\tilde{y}(t)$ given in \eqref{eq:err_dyn2} are considered in the performance output $z^\mrm{OE}(t)$. This improves the feasibility of the method as the choice of $c_2>0$ is no longer constrained by the size of $Q^\top Q$. Following this idea, we have the next result.
	\begin{corollary}
		\label{co:tuning2}
		Consider the system \eqref{eq:err_dyn_dis} together with Assumption~\ref{ass:1}, input $\delta(t)$ and output $z^\mrm{OE}(t)$ given in \eqref{eq:z_OE}. Choose $\alpha>0$ as
		\begin{equation}
			\label{eq:alpha_tuning2}
			\underset{\alpha}{\mrm{arg \ max}} \ \frac{\Ml(\alpha)}{\Mu(\alpha)},
		\end{equation}
		with $\Ml$, and $\Mu$ defined in \eqref{eq:iota_3}. Fix the constant $c_2>0$ and solve the following \ac{sdp} \eqref{eq:sdp2} in the scalar decision variables ($c_1, \gamma_1, \gamma_2, \gamma$):
		\begin{equation}
			\label{eq:sdp2}
			\begin{aligned}
				&\mrm{minimize} \ \ \gamma&\\
				&\mrm{subject \ to:}	 	& \gamma >0, \ \gamma_1 - \gamma_2 \ge 0, \\
				&				&\begin{bmatrix}
					\Phi_{11}^\mrm{OE} & \Phi_{12}\\ \Phi_{12}^\top & -\gamma_2 I_{nN}
				\end{bmatrix} &\le0,
			\end{aligned}
		\end{equation}
		where $\Phi_{11}^\mrm{OE}$ and $\Phi_{12}$ are defined in \eqref{eq:Phi}. If the \ac{sdp} \eqref{eq:sdp2} is feasible, choose $\bar{\Gamma}$, such that
		\begin{equation}
			\label{eq:LMI_Gamma2}
			\sqrt{\gamma_1} I_{nN } \ge \bar{\Gamma} \ge  	\sqrt{\gamma_2} I_{nN }.
		\end{equation}
		In such case, the system \eqref{eq:err_dyn_dis} with output \eqref{eq:z_OE} has an $\mathcal{L}_2$-gain of at most $\sqrt{\gamma}$.
	\end{corollary}

	In view of the results in Corollaries \ref{corollary:convergence_bounds} to \ref{co:tuning2}, we propose a two-step tuning process. First, $\alpha$ satisfying \eqref{eq:alpha_tuning2} is chosen in a way that the convergence time (see \eqref{eq:trajecty_bound}) and the \ac{iss}-gain (see \eqref{eq:ISS_gain}) are minimized. This requires an initial guess of $\bar{\Gamma}_0$, since the convergence time and the \ac{iss}-gain depend on $\bar{\Gamma}$. Then, $\bar{\Gamma}$ is updated by solving the \ac{sdp} \eqref{eq:sdp2}. Hence, guaranteeing an $\mathcal{L}_2$-gain of at most $\sqrt{\gamma}$. 

	In a practical setting, the excitation parameters $\Mu$, $\Ml$, $r_4$, and $T$ are typically not known exactly but rather upper and lower bounds for them, e.g., obtained by analyzing previous measurement data. The following proposition is a direct extension of \cite[Proposition 1]{rueda-escobedo_l2-gain_2022} to the tuning method presented in this work. It shows that exact knowledge of the excitation parameters is not required, but upper and lower bounds are sufficient to guarantee robust performance. 
	\begin{proposition}
		\label{pro:tuning}
		Assume that \eqref{eq:LMI_tuning_th1} is satisfied with constants $\gamma^*>0$, $\gamma_1^* \ge\gamma_2^*>0$, $\alpha^*>0$, $c_1^*>0$, $c_2^*>0 $, $\Mu^*\ge\Ml^*>0$, $r_4^*>0$, and $T^*>0$. Then, \eqref{eq:LMI_tuning_th1} remains valid for any $\Mu$, $\Ml$, $r_4$, and $T^*>0$ that satisfy 	
		\begin{equation}
			\label{eq:bounds_excitation}
			\Mu^*\ge\Mu\ge \Ml\ge \Ml^*>0, \ \ r_4^*\ge r_4 >0, \ \  T^*\ge T>0.
		\end{equation}
	\end{proposition}
	It is straightforward to extend this proposition to the results of Corollary~\ref{co:tuning} and~\ref{co:tuning2}. Hence, these additions are not stated here explicitly. Moreover, as as seen from Lemmas~\ref{lemma:obs_gram_OI} and \ref{lemma:obs_gram}, $\Mu$ and $\Ml$ depend on the choice of gains $\alpha$ and $\bar{\Gamma}$. For the application of Theorem~\ref{th:tuning} as well as Corollaries~\ref{co:tuning} and \ref{co:tuning2}, an admissible range of possible gains should be decided on a priori. Hence, $\Mu^*$ and $\Ml^*$ can be calculated using the limit cases of the selected range.
	\section{Application example}
	\label{sec:acc_ex}
	To illustrate our results, we consider six identical mass-spring-damper systems (cf. \cite{chen_distributed_2014}). The $i$th system model is given by 
	\begin{equation}
		\begin{split}		
			\label{eq:sys_mass_spring}
			\dot{\xi}_{1,i}(t) &= \xi_{2,i}(t), \\
			\dot\xi_{2,i}(t) &= \frac{1}{k_1} \left(u_i(t) -k_2\xi_{1,i}(t) - k_3(t)\xi_{2,i}(t)\right),
	\end{split}\end{equation}
	with $\dot\xi_{2,i}(t)$, $\xi_{2,i}(t)$, $\xi_{1,i}(t)$, and $u_i(t)$ denoting the known local acceleration, velocity, displacement, and input, respectively. The unknown and constant global spring coefficient and mass are denoted by $k_2$ and $k_1$, respectively. The damper coefficient $k_3(t)$ is considered slowly time-varying and modeled as $\dot{k}_3(t) = d_1\sin(0.5t)$, with the amplitude $d_1$. The inputs are defined as $u_1(t) = \sin(t)$, $u_2(t)= 2\cos(0.5t)$, $u_3(t) = 3\sin(3t)$, $u_4(t) =3\cos(2t)$, $u_5(t) =\sin(t)+0.5\cos(t)$, and $u_6(t) =2\sin(3t)+\cos(0.4t)$. 
	\begin{figure}
		\centering
		\includegraphics[width=0.2\textwidth]{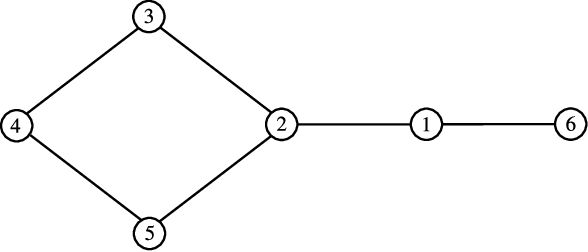}
		\caption{Communication structure among the six agents.}
		\label{fig:com_stru}
	\end{figure}
	The communication structure among the agents is depicted in Figure~\ref{fig:com_stru}.
	
\begin{figure}
		\begin{subfigure}[t]{0.24\textwidth}
			\includegraphics[width=\textwidth]{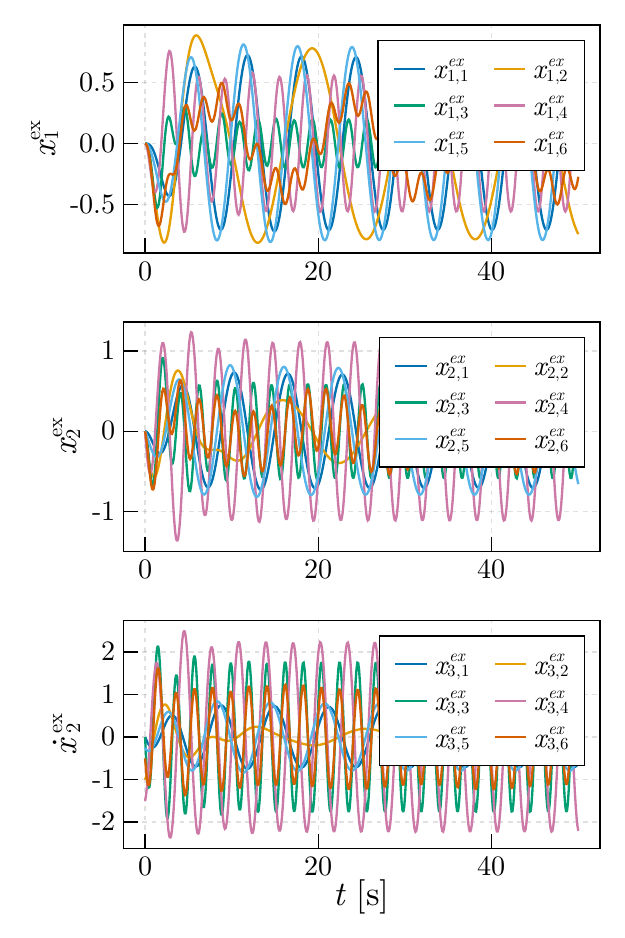}
			\caption{\tabular[t]{@{}l@{}}Constant parameters: \\ $d_1=0$\endtabular}
		\end{subfigure}\hfill
		\begin{subfigure}[t]{0.24\textwidth}
			\includegraphics[width=\textwidth]{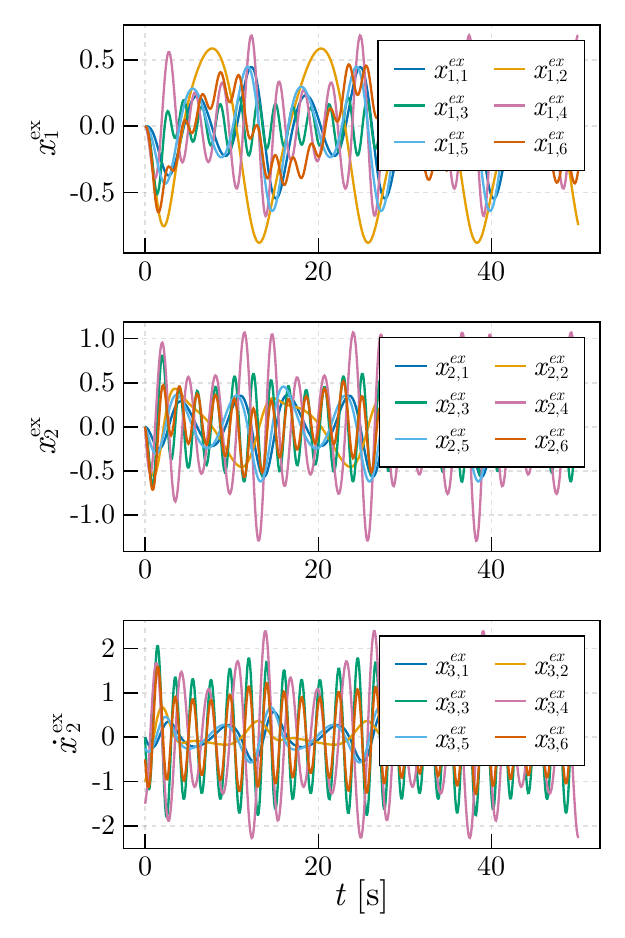}
			\caption{\tabular[t]{@{}l@{}}Time-varying parameters:\\ $d_1=1$ \endtabular}
		\end{subfigure}
		\caption{Simulation results of the system \eqref{eq:sys_mass_spring}. }
		\label{fig:sys_sim}
	\end{figure}

	We simulate the system \eqref{eq:sys_mass_spring} for 50 seconds using the programming language Julia and the package \mbox{DifferentialEquations.jl} \cite{rackauckas_differentialequationsjl_2017}. The results for $d_1=0$ and $d_1=1$ are shown in Figure~\ref{fig:sys_sim}.
	
	We facilitate the application of the distributed parameter estimator \eqref{eq:C+I_alg} by recasting the system \eqref{eq:sys_mass_spring} as a \ac{lre} \eqref{eq:local_LRE} with 
	\begin{equation}
		\label{eq:lre_example}
		\begin{split}
			y_i^\mrm{ex}(t) &= \dot\xi_{2,i}(t), \quad C_i^\mrm{ex}(t) = \begin{bmatrix}u_i(t) & -\xi_{1,i}(t) &- \xi_{2,i}(t)
			\end{bmatrix}, \\\theta^\mrm{ex}(t) &= \begin{bmatrix}1/k_1 & k_2/k_1 & k_3(t)/k_1
			\end{bmatrix}^\top.
		\end{split}
	\end{equation}
	Furthermore, we consider the disturbance $\delta(t)\in \mathbb{R}^3$
	\begin{equation}
		\label{eq:disturbance}
		\delta(t) = \begin{bmatrix}
			d_1 \sin(0.5t) &d_2 \sin(50t) & d_3 \sin(50t)
		\end{bmatrix}^\top,
	\end{equation}
	where the first component models the parameter variation, the second one local measurement noise, and the third one disturbances in the communication channels. This is reflected in the disturbed error dynamics \eqref{eq:err_dyn_dis} by defining 
	\begin{equation*}
		\Db1 = \bm{1}_6\otimes \begin{bmatrix}
			0 & 0 & 0\\
			0 & 0 & 0\\
			0 & 0 & 1
		\end{bmatrix}, \quad \Db2 = \begin{bmatrix}
			\bm{0}_6 & \bm{1}_6 & \bm{0}_6 \\\bm{0}_{18} & \bm{0}_{18} & \bm{1}_{18} 
		\end{bmatrix}.
	\end{equation*}
	To employ the tuning method as specified in Corollary~\ref{co:tuning2}, we define the performance output $z \in \mathbb{R}^5$ as in \eqref{eq:z_OE} and choose
	\begin{equation*}
		Q^\mrm{OE} = \begin{bmatrix}
			\bm{1}_6 & \bm{0}_6 & 0_{6 \times 3} \\\bm{0}_{18} & \bm{1}_{18} & 0_{{18} \times 3} 
		\end{bmatrix}^\top ,\quad  W^\mrm{OE} = \begin{bmatrix}0_{2 \times 3} \\ I_3
		\end{bmatrix}.
	\end{equation*}
	For applying Corollary~\ref{co:tuning2}, lower and upper bounds on the observability Gramian \eqref{eq:gram} are required (see Proposition~\ref{pro:tuning}). Although the existence of these bounds is proven in Lemma~\ref{lemma:obs_gram}, the exact values provided are conservative. This makes it more difficult to solve the optimization problem formulated in Corollary~\ref{co:tuning2}. We obtain less conservative values by approximating these bounds by numerical calculation of the observability Gramian \eqref{eq:gram}. For this, we differentiate the observability Gramian w.r.t. time to obtain its differential equation 
	\begin{equation}
		\begin{split}
			\label{eq:gram_ode}
			\dot{M}(t) &=\la^\top(t) \la(t) \bar{\Gamma} M(t) + M(t) \bar{\Gamma}\la^\top(t) \la(t) \\& \ \ \ + \la^\top(t) \la(t) \ \ , \quad  \qquad M(t_0) = 0.
		\end{split}
	\end{equation}
	As $\la(t)$ depends on the simulation results of the system \eqref{eq:sys_mass_spring}, the oriented incidence matrix $D(t)$, and the gain $\alpha$ (see \eqref{eq:y_hat}), we can calculate $M(t)$ for different initial times $t_0$, integration intervals $T$, and gains $\bar{\Gamma}$ and $\alpha$. The lower and upper bounds on the resulting Gramian $M(t=T)$ are obtained by calculating its smallest and largest eigenvalue, respectively. 
	
	After simulating the mass-spring-damper systems \eqref{eq:sys_mass_spring}, we numerically approximate these bounds by solving \eqref{eq:gram_ode} with a fixed integration interval $T=0.01$ for 1000 different initial times. The results are searched for the worst case. We estimate the lower bound $\Ml$ using the smallest values of the gains that we expect to be reasonable for the system setup. For an estimate of the upper bound $\Mu$, we use the largest expected values of the gains. 
	
	After obtaining $\Ml$ and $\Mu$ as explained above, we derive gains $\bar{\Gamma}$ and $\alpha$, which minimize the $\mathcal{L}_2$-gain by following Corollary~\ref{co:tuning2}. We solve the optimization problem utilizing the Julia package JuMP \cite{lubin_jump_2023} and solve the \ac{sdp} with MOSEK \cite{mosek_aps_mosek_2023}. The optimal gains found by this procedure are $\bar{\Gamma}^*= 1.93I_{nN} $ and $\alpha^* = 1.05$. 
	
\begin{figure}
		\begin{subfigure}[t]{0.24\textwidth}
			\includegraphics[width=\textwidth]{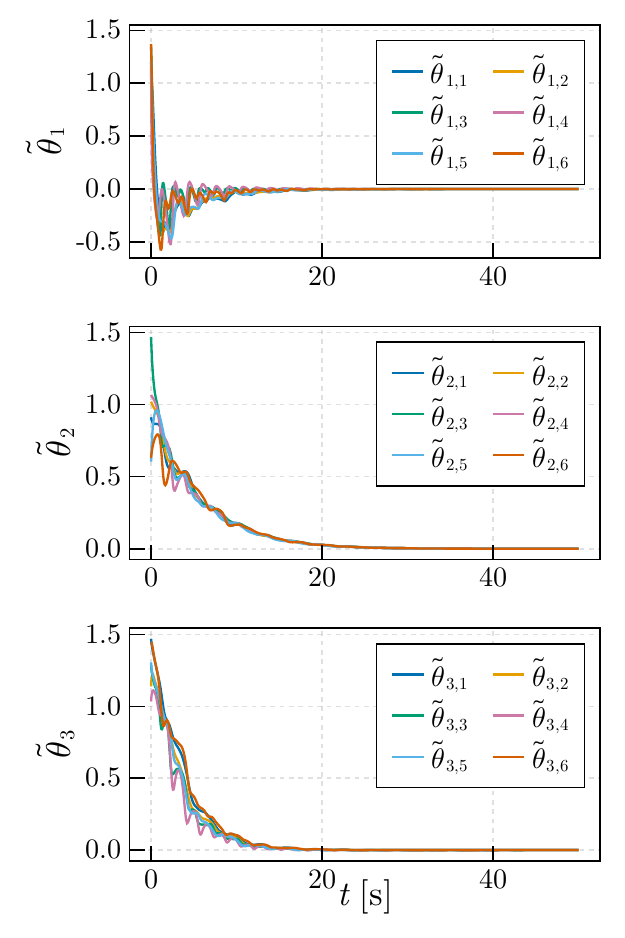}
			\caption{\tabular[t]{@{}l@{}}Nominal scenario: \\ $d_1=d_2=d_3=0$\endtabular }
		\end{subfigure}\hfill
		\begin{subfigure}[t]{0.24\textwidth}
			\includegraphics[width=\textwidth]{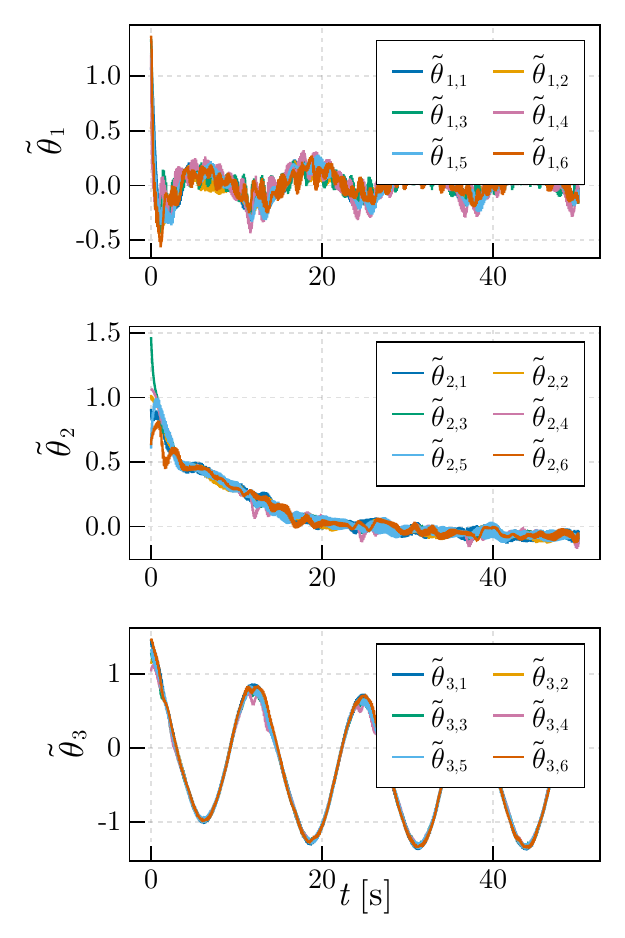}
			\caption{\tabular[t]{@{}l@{}}Scenario 2: \\ $d_1= 0.5$, $d_2=1$, $d_3=0.5$\endtabular }
		\end{subfigure}
		\caption{Estimation errors of the six agents. }
		\label{fig:CI_est}
	\end{figure}

	To investigate the performance of the algorithm using the optimized gains, we define five scenarios with varying disturbances acting on the system, modeled by different amplitudes $d_1$, $d_2$, and $d_3$ of the disturbance components \eqref{eq:disturbance}. The considered scenarios are detailed in Table~\ref{tab:scenarios}. 
	\begin{table}
		\centering
		\begin{tabular}{c|c|c|c|l}
			\cmidrule{1-5}
			Scenario & $d_1$& $d_2$ &$d_3$ & Description \\  \cmidrule{1-5}\morecmidrules\cmidrule{1-5}	
			\multirow{2}{*}{1}  & \multirow{2}{*}{$0$} & \multirow{2}{*}{$1$} & \multirow{2}{*}{$\frac{1}{2}$} & No parameter variation, high noise level, \\
			&  & & & high communication disturbances \\\cmidrule{1-5}
			\multirow{2}{*}{2}  & \multirow{2}{*}{$\frac{1}{2}$} & \multirow{2}{*}{$1$} & \multirow{2}{*}{$\frac{1}{2}$} & Slow parameter variation, high noise level, \\
			&  & & & high communication disturbances \\\cmidrule{1-5}
			\multirow{2}{*}{3}  & \multirow{2}{*}{$2$} & \multirow{2}{*}{$\frac{1}{4}$} & \multirow{2}{*}{$\frac{1}{8}$} & Fast parameter variation, low noise level, \\
			&  & & & low communication disturbances \\\cmidrule{1-5}
			\multirow{2}{*}{4}  & \multirow{2}{*}{$2$} & \multirow{2}{*}{$0$} & \multirow{2}{*}{$0$} & Fast parameter variation, no noise level, \\
			&  & & & no communication disturbances \\\cmidrule{1-5}
			\multirow{2}{*}{5}  & \multirow{2}{*}{$2$} & \multirow{2}{*}{$1$} & \multirow{2}{*}{$\frac{1}{2}$} & Fast parameter variation, high noise level, \\
			&  & & & high communication disturbances \\\cmidrule{1-5}
			\cmidrule{1-5}
		\end{tabular}
		\caption{Considered scenarios with the corresponding disturbance amplitudes.}
		\label{tab:scenarios}
	\end{table}
	We illustrate the trajectories of the estimation errors for a nominal scenario without perturbations and Scenario 2 according to Table~\ref{tab:scenarios} in Figure~\ref{fig:CI_est}. The initial conditions of the parameter estimators are chosen randomly. The estimation errors of the six agents shown in Figure~\ref{fig:CI_est} are defined as 
	$$\tilde{\theta}^\mrm{ex}_{i,j}(t) = \hat{\theta}^\mrm{ex}_{i,j}(t)- \theta^\mrm{ex}_{i,j}(t),$$ where \mbox{$i = \{1,2,3\}$} denotes the considered parameter and \mbox{$j=\{1,2,...,6\}$} denotes the agent estimating the parameter. 
	We obtain $\hat{\theta}^\mrm{ex}_{i,j}(t)$ by solving \eqref{eq:C+I_alg} and $\theta^\mrm{ex}_{i,j}(t)$ is defined in \eqref{eq:lre_example}. 
	
\begin{figure}
		\centering
		\includegraphics[width=0.5\textwidth]{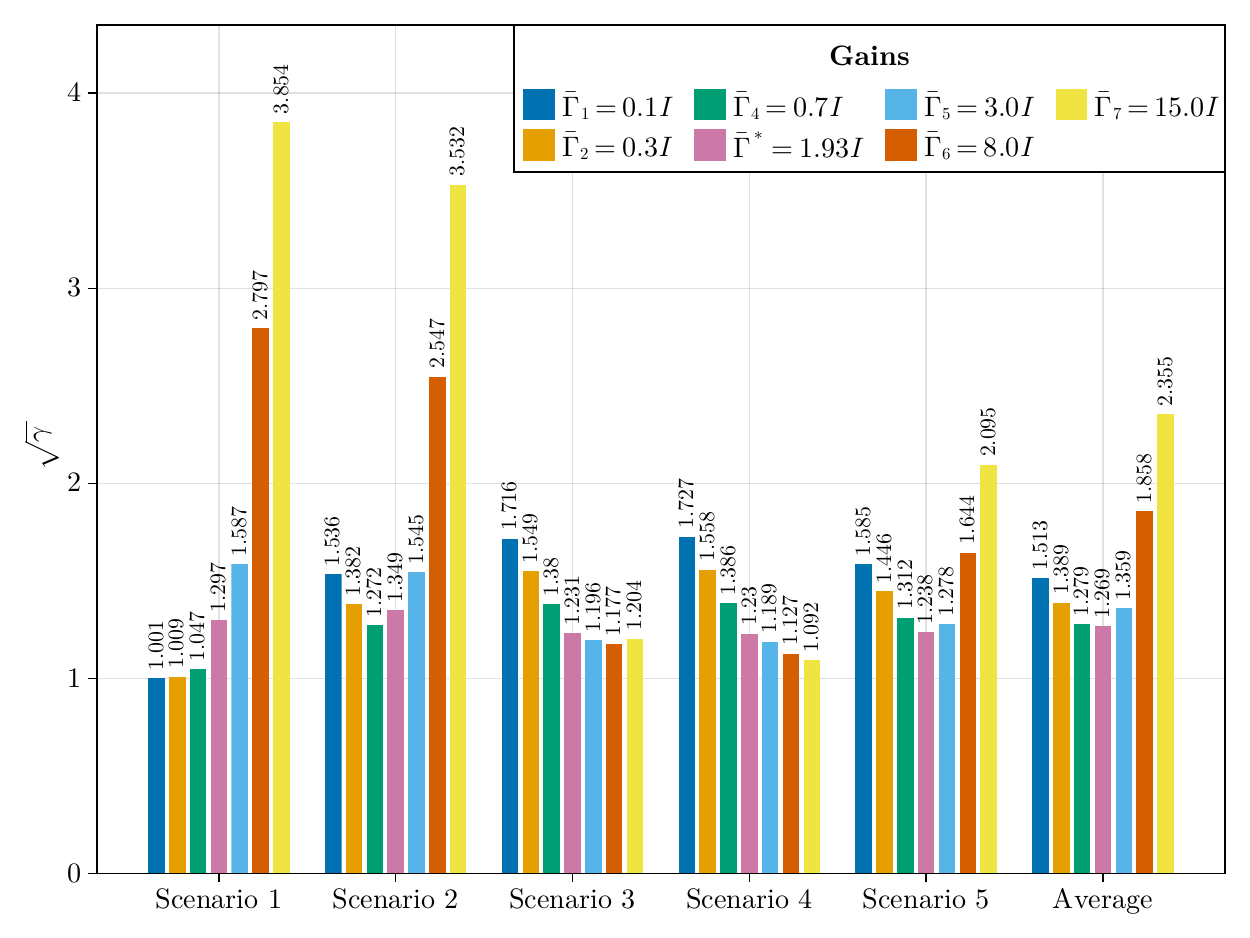}
		\caption{The calculated values of the metric $\sqrt{\gamma}$ as defined in \eqref{eq:metric} for the five scenarios given in Table~\ref{tab:scenarios} using different gains, as well as the average value of $\sqrt{\gamma}$ for each gain over the five scenarios.}
		\label{fig_histogram}
	\end{figure}

	The results of all scenarios are evaluated using the $\mathcal{L}_2$-gain as a metric 
	\begin{equation}
		\label{eq:metric}
		\sqrt{\gamma} = \frac{\sqrt{\int_{0}^{50}\left(z^\mrm{OE}\right)^\top(s)z^\mrm{OE}(s)\mrm{d}s}}{\sqrt{\int_{0}^{50}\delta^\top(s)\delta(s)\mrm{d}s}}.
	\end{equation}
	For each of the five scenarios, we simulate the disturbed estimation error dynamics \eqref{eq:err_dyn_dis} with zero initial conditions and the optimized gains as well as six additional gains chosen smaller and larger than the optimal ones. The values of the metric $\sqrt{\gamma}$ as defined in \eqref{eq:metric} for the different scenarios and gains, as well as the average value of $\sqrt{\gamma}$ for each gain over the five scenarios are shown in Figure~\ref{fig_histogram}.
	
	As seen from Figure~\ref{fig_histogram}, the optimized gains $\alpha^*$ and $\bar{\Gamma}^*$ provide a good trade-off between fast convergence speed, measurement noise, and communication disturbance attenuation and perform well in all scenarios. The optimized gains show the smallest average value of $\sqrt{\gamma}$ over all five scenarios of all the gains. In Scenario 1, with high noise levels and high communication disturbances but no parameter variations, low gains perform better. However, once parameter variations are considered, the performance deteriorates. Without noise and communication disturbances (Scenario 4), higher gains improve the performance, whereas they perform poorly in the presence of noise and communication disturbances. In Scenario 5, with fast parameter variations, high noise levels, and high communication disturbances, the optimized gains show the smallest value of $\sqrt{\gamma}$, demonstrating the effectiveness of the proposed tuning method.

	\section{Conclusions}
	\label{sec:conclusions}
	In this paper, we recast the error dynamics of a \ac{c+i}-based distributed parameter estimator to reflect the error dynamics induced by the classical gradient descent algorithm. Proving the \ac{uco} of this system paved the way to obtain a strong \ac{lf} by using a scaled version of the observability Gramian as a strictifying term. This strong \ac{lf} was utilized to derive convergence bounds for the error dynamics and the \ac{iss}-gain w.r.t. an affine disturbance. Also, it allowed us to provide \ac{lmi}-based tools to tune the algorithm's gains such that a guaranteed upper bound on the $\mathcal{L}_2$-gain with respect to parameter variations, measurement noise, and disturbances in the communication channels is achieved. Finally, an application example demonstrated the advantages of the derived method. It was shown that the optimized gains provide a good trade-off between fast convergence speed and the attenuation of measurement noise and communication disturbances. 
	The derived tuning criteria greatly benefit the practical application of \ac{c+i}-based distributed parameter estimators, and it is ongoing work to apply this method to power system applications.
	
	\section{Proof of the main results}
	\label{sec:proofs}
	\subsection{Proof of Lemma \ref{lemma:obs_gram_OI}}
	\begin{proof}[\unskip\nopunct]
		For establishing the result of Lemma~\ref{lemma:obs_gram_OI}, we need to derive a lower and upper bound on  the observability Gramian $M^\mrm{OI}(t,t-T)$ as defined in \eqref{eq:gram_OI}, which can be expanded as 
		\begin{gather}
			M^\mrm{OI}(t,t-T)= \bar{H}(t)+\bar{K}(t),\label{eq:F_bar_bound}\\
			\bar{H}(t) = \int_{t-T}^{t}\bar{C}^\top(s)\bar{C}(s) \mrm{d}s, \qquad \bar{K}(t) = \int_{t-T}^{t}\alpha\bar{L}(s)\mrm{d}s.\nonumber
		\end{gather}
		The existence of a lower bound is proven in \cite{chen_distributed_2014}, and a lower bound for time-invariant graphs is provided in \cite{papusha_collaborative_2014}. Building on these results, we extend the result to time-varying graphs.
		Following Assumption~\ref{ass:1}, $\int_{t-T}^{t}L(s)\mathrm{d}s$ has only one zero eigenvalue $\lambda_1$ with the corresponding constant unit eigenvector $\nu^\mrm{sp}_1 = \frac{1}{\sqrt{n}}\bm{1}_n$. From the spectrum property of the Kronecker product \cite[Chaper 8.2]{bullo_lectures_2020}, it follows that $\int_{t-T}^{t}\bar{L}(s)\mrm{d}s$ has $N$ zero eigenvalues $\lambda_1, ... \lambda_N$, whose $N$ constant eigenvectors can be written as $\nu_i = \frac{1}{\sqrt{n}}\bm{1}_n\otimes e_i$, where $e_i\in\mathbb{R}^N$ denotes a standard Cartesian unit vector. The remaining $nN-N$ time-varying unit eigenvectors $\mu_i(t)$ have associated positive time-varying eigenvalues, which can be lower bounded by the constant  $\underline{\lambda} \le \lambda_i(t) , \  \forall i \ge N+1$. Consequently, we can express any unit vector $w$ as a linear combination of these eigenvectors as 
		\begin{equation}
			\label{eq:linear_combination_EV_L}
			w = \sum_{i=1}^{N} a_i(t) \nu_i + \sum_{i=N+1}^{nN} b_i(t) \mu_i(t),
		\end{equation} 
		with $a(t) \coloneqq [a_1(t), ..., a_N(t)]^\top$, $b(t) = [b_{N+1}(t), ... , b_{nN}]^\top$, and $\|a(t)\|^2+  \|b(t)\|^2 = 1 $.
		\\Defining $V \coloneqq [\nu_1, ..., \nu_N]$, we can express the first term of \eqref{eq:F_bar_bound} as
		\begin{equation*}
			\begin{split}
				&w^\top \bar{H}(t) w = a^\top(t) V^\top \bar{H}(t) V a(t) \\ &  + 2 \sum_{i=0}^{N} a_i(t) \nu_i^\top \bar{H}(t) \sum_{j=N+1}^{nN} b_j(t) \mu_j(t)\\
				& + \underbrace{\sum_{i=N+1}^{nN} b_i(t) \mu_i(t)^\top\bar{H}(t) \sum_{j=N+1}^{nN} b_j(t) \mu_j(t)}_{\ge 0},
			\end{split}
		\end{equation*}
		where the third term is non negative since $\bar{C}^\top(t)\bar{C}(t)\ge 0$. The first term can be bounded employing the \ac{cpe} \eqref{eq:cPE} and defining $ H_i(t) \coloneqq \int_{t-T}^{t}C_i^\top(s)C_i(s)\mrm{d}s$ as 
		\begin{equation*}
			a^\top(t)	V^\top \bar{H}(t) Va(t)  \ge \frac{\cPEl}{n} \|(t)\|^2,
		\end{equation*}
		where the derivation is provided in \eqref{eq:VHV}.
		\begin{figure*}[!t]
			\scriptsize
			\begin{equation}
				\label{eq:VHV}
				\begin{split}
					& a^\top(t)	V^\top \bar{H}(t) Va(t) = \frac{1}{n}a^\top(t)\begin{bmatrix}
						e_1^\top & ... & e_1^\top \\
						\vdots   & \ddots & \vdots \\
						e_N^\top & ... & e_N^\top 
					\end{bmatrix}
					\begin{bmatrix}
						H_1(t) & 0 &... & 0 \\
						0  &	H_2(t)   &... & 0 \\
						\vdots   & \vdots   & \ddots & \vdots \\
						0  &	0  &... & H_n(t) 
					\end{bmatrix}
					\begin{bmatrix}
						e_1 & ... & e_N \\
						\vdots   & \ddots & \vdots \\
						e_1 & ... & e_N 
					\end{bmatrix}a(t)
					\\& =\frac{1}{n} a^\top(t)\begin{bmatrix}
						e_1^\top \sum_{i=1}^{n}H_i(t)e_1 & e_1^\top \sum_{i=1}^{n}H_i(t)e_2 &... & e_1^\top \sum_{i=1}^{n}H_i(t)e_N \\
						e_2^\top \sum_{i=1}^{n}H_i(t)e_1 & e_2^\top \sum_{i=1}^{n}H_i(t)e_2 &... & e_2^\top \sum_{i=1}^{n}H_i(t)e_N \\
						\vdots 	&\vdots   & \ddots & \vdots \\
						e_N^\top \sum_{i=1}^{n}H_i(t)e_1 & e_N^\top \sum_{i=1}^{n}H_i(t)e_2 &... & e_N^\top \sum_{i=1}^{n}H_i(t)e_N 
					\end{bmatrix}a(t)
					= \frac{1}{n}a^\top(t)	\int_{t-T}^{t}\sum_{i=1}^{n}{C_i}^\top(s){C_i}(s)\mathrm{d}sa(t) \ge \frac{\cPEl}{n} \|a(t)\|^2.
				\end{split}
			\end{equation}
			\begin{equation}
				\label{eq:VHU}
				\begin{split}
					&	\left|2 \sum_{i=0}^{N} a_i(t) \nu_i^\top \bar{H}(t) \sum_{j=N+1}^{nN} b_j(t) \mu_j(t)\right| \le 2 \left\|a^\top(t)V^\top\right\| \left\|\bar{H}(t)\right \| \left\|U(t)b(t)\right \| \le  2 	\left\|a(t)\right \| \underbrace{\left\|V\right \|}_{= 1} 
					\underbrace{\left\|U(t)\right \|}_{= 1} \left\|b(t)\right \|	\left\|\begin{bmatrix}
						H_1(t) & 0 &... & 0 \\
						0  &	H_2(t)   &... & 0 \\
						\vdots   & \vdots   & \ddots & \vdots \\
						0  &	0  &... & H_n(t) 
					\end{bmatrix} \right\|  \\
					&\le  2 	\left\|a(t)\right \| \left\|b(t)\right\| \underset{\forall t,i}{\mrm{max}} \left\|	H_i(t)	\right\| 
					=  2 	\left\|a(t)\right \| \left\|b(t)\right\| \underset{\forall t,i}{\mrm{max}} \left\|	\int_{t-T}^{t}C_i^\top(s)C_i(s)\mrm{d}s	\right\| 	\le  2 	\left\|a(t)\right \| \left\|b(t)\right\| \underset{\forall t,i}{\mrm{max}} 	\int_{t-T}^{t}\left\|C_i^\top(s)C_i(s)\right\|\mrm{d}s	 \\ &\le  2 r_2 T	\left\|a(t)\right \|  \left\|b(t)\right \|	 = 2 r_2 T \sqrt{ \left\| a(t)\right\|^2 (1-\left\|a(t)\right \|^2)}.
				\end{split}
			\end{equation}
			\hrule 
		\end{figure*}
		The second, indefinite term can be bounded by defining the orthogonal matrix $U(t) \coloneqq [\mu_{N+1}(t), ..., \mu_{nN}(t)]$, recalling Assumption~\ref{ass:1}, and by using the Cauchy-Schwarz inequality \cite[Equation B7]{horn_matrix_2012} 
		\begin{equation*}
			\begin{split}
				&	\left|2a^\top(t) V^\top \bar{H}(t) U(t)b(t)\right| = 2 r_2 T \sqrt{ \left\| a(t)\right\|^2 (1-\left\|a(t)\right \|^2)}.
			\end{split}
		\end{equation*}
		The derivation is provided in \eqref{eq:VHU}. 
		Hence, we have 
		\begin{equation*}
			w^\top \bar{H}(t) w \ge \frac{\cPEl}{n} \|a(t)\|^2 -    2 r_2 T  \sqrt{ \left\| a(t)\right\|^2 (1-\left\|a(t)\right \|^2)}.
		\end{equation*}
		For the second term of \eqref{eq:F_bar_bound}, we derive the following bound 
		\begin{equation*}
			\begin{split}
				&w^\top \bar{K}(t) w = \alpha\sum_{i=0}^{N} a_i(t) \nu_i^\top  \underbrace{\int_{t-T}^{t}\bar{L}(s)\mrm{d}s \sum_{j=0}^{N} a_j \nu_j}_{=\bm{0_{nN}}}  \\&
				+ 2 \alpha \sum_{i=N+1}^{nN} b_i(t) \mu_i(t)^\top\underbrace{\int_{t-T}^{t}\bar{L}(s)\mrm{d}s\sum_{j=1}^{N} a_j(t) \nu_j}_{=\bm{0_{nN}}}  \\
				& +		\alpha\sum_{i=N+1}^{nN} b_i(t) \mu_i(t)^\top \int_{t-T}^{t}\bar{L}(s)\mrm{d}s\sum_{j=N+1}^{nN} b_j(t) \mu_j(t)\\ &		
				\\		&= + \alpha\sum_{i=N+1}^{nN} \sum_{j=N+1}^{nN}  b_i(t)  b_j(t) \lambda_j(t)\underbrace{\mu_i(t)^\top  \mu_j(t)}_{\substack{\ \  =1,\ \mrm{if} \ i=j,\\  =0, \ \mrm{else.}}} \\
				&=  \alpha\sum_{i=N+1}^{nN} b_i^2(t)  \lambda_i(t) \ge \alpha  \underline{\lambda} \|b(t)\|^2 = \alpha  \underline{\lambda} ( 1-\|a(t)\|^2 ),
			\end{split}
		\end{equation*}
		and consequently, we establish the following strictly positive lower bound 
		\begin{equation*}
			\begin{split}
				&	M^\mrm{OI}(t,t-T) \ge  \underset{\forall \left\| a(t)\right\|^2 \le 1}{\mrm{min}} \mrm{max} \bigg\{ \alpha  \underline{\lambda} ( 1-\|a(t)\|^2 ), \\ & \ \ \ \frac{\cPEl}{n} \|a(t)\|^2 -    2 r_2 T  \sqrt{ \left\| a(t)\right\|^2 (1-\left\|a(t)\right \|^2)} \bigg\}  \eqqcolon \MOIl,
			\end{split}
		\end{equation*}	
		as $\bar{H}(t)$ and $\bar{K}(t)$ are positive semi-definite.\\
		To establish a upper bound of $M^\mrm{OI}(t,t-T)$, we recall Assumption~\ref{ass:1} and bound the block-diagonal matrix $\bar{C}(t)^\top \bar{C}(t)$ with blocks ${C}_i(t)^\top {C}_i(t)$ as 
		\begin{equation*}
			\|\bar{C}(t)^\top \bar{C}(t)\| \le   \underset{\forall t, i}{\mrm{max}} \ \lambda_\mrm{max}(C_i(t)^\top C_i(t)) = r_2.
		\end{equation*}
		Moreover, from the spectrum property of the Kronecker product, we have that
		$$\|\bar{L}(t)\| = \|L(t)\| \le \underset{\forall t}{\mrm{max}} \ \lambda_\mrm{max}(L(t)) = r_3.$$
		Hence, we have that 
		\begin{equation*}
			M^\mrm{OI}(t,t-T) \le \int_{t-T}^{t} r_4 \mrm{d}s = Tr_4 \eqqcolon \MOIu,
		\end{equation*}
		which concludes the proof.
	\end{proof}
	\subsection{Proof of Lemma \ref{lemma:obs_gram}}
	\begin{proof}[\unskip\nopunct]
		We first establish a bound on the output matrix $\la(t)$ as
		\begin{equation*}
			\begin{split}
				\|\la(t)\| &= \|\la^\top(t)\| = \|\la^\top(t) \la(t) \|^\frac{1}{2} \\ &=\|\bar{C}^\top(t) \bar{C}(t) +\alpha \bar{L}(t)\|^\frac{1}{2} \le \sqrt{r_4}.
			\end{split}
		\end{equation*}
		Then, we follow the proof of \cite{zhang_observability_2015} until Equation 23. As the system matrix of \eqref{eq:err_dyn2_OI} is the zero matrix $A(t) = 0$, it has the bound $\|A(t)\| \le 0$ (denoted as $\gamma$ in \cite{zhang_observability_2015}), making \cite[Equation 24]{zhang_observability_2015} ill-conditioned. As the state-transition matrix of \eqref{eq:err_dyn2_OI} is the identity matrix, we rewrite \cite[Equation 23]{zhang_observability_2015} as
		\begin{equation*}
			\int_{s}^{t} \|I_{nN}\|^2 \mrm{d}p = t-s.
		\end{equation*} 
		Then, following the procedure of Zhang et al. and denoting by $w$ an arbitrary unit vector, \cite[Equation 25]{zhang_observability_2015} becomes 
		
		\begin{equation*}
			\begin{split}
				&\int_{t-T}^{t}  \left\| \la(s)\int_{s}^{t}  I_{nN}\bar{\Gamma}\la^\top(p) \la(p) \bar{\Phi}(p,t)w  \mathrm{d}p \right\|^2  \mathrm{d}s 
				\\
				\\
				& \le r_1^2r_4^2 \int_{t-T}^{t}   \int_{s}^{t}  \mathrm{d}p\int_{s}^{t} \| \la(p) \bar{\Phi}(p,t)w\|^2  \mathrm{d}p \mathrm{d}s
				\\ & = \frac{1}{2}r_1^2r_4^2T^2  \int_{s}^{t} \| \la(p) \bar{\Phi}(p,t)w\|^2  \mathrm{d}p.
			\end{split} 
		\end{equation*}
		Hence, and differing from \cite{zhang_observability_2015}, we define 
		$$\varphi_1 \coloneqq\frac{1}{2}r_1^2r_4^2T^2,$$
		which is denoted as $\varphi$ in \cite{zhang_observability_2015}. The rest of the proof remains equivalent to \cite{zhang_observability_2015}.
	\end{proof}

	\subsection{Proof of Theorem \ref{th:stong_lya}}
	\begin{proof}[\unskip\nopunct]
		The first term of the \ac{lf} candidate $V(\tilde{\Th}, t)$ \eqref{eq:lya_fun} can be bounded as 
		\begin{equation*}
			\begin{split}
				\frac{T}{2}\lambda_\mrm{max}(\bar{\Gamma}^{-1}) \|\tilde{\Th}(t)\|^2 &\ge	\tilde{\Th}^\top(t) \frac{T}{2} \bar{\Gamma}^{-1} \tilde{\Th}(t) \\&\ge \frac{T}{2} \lambda_\mrm{min}(\bar{\Gamma}^{-1}) \|\tilde{\Th}(t)\|^2.
			\end{split}	
		\end{equation*}
		For ease of notation, we define 
		$$f(s,t)\coloneqq  \bar{\Phi}^\top(s,t)\la^\top(s)\la(s)\bar{\Phi}(s,t).$$ As $\Mu I_{nN}\ge\int_{t-T}^{t}f(s,t)\mathrm{d}s $ (see Lemma~\ref{lemma:obs_gram}), and ${0\le s-t+T \le T}$ in the time-interval $s \in [t-T, T]$, we can bound the second term of $V(\tilde{\Th}, t)$ as 
		\begin{equation*}
			\begin{split}
				T\Mu\|\tilde{\Th}(t)\|^2 
				\ge\tilde{\Th}^\top(t)\int_{t-T}^{t}(s-t+T)f(s,t) \mathrm{d}s	\tilde{\Th}(t) \ge	0 .
			\end{split}
		\end{equation*} 
		Thus, the \ac{lf} candidate \eqref{eq:lya_fun} is positive definite with the following bounds
		\begin{equation}
			\begin{split}
				\label{eq:V_bounds}
				&	\left( \frac{T}{2}\lambda_\mrm{max}(\bar{\Gamma}^{-1}) + T\Mu\right)\|\tilde{\Th}(t)\|^2 \ge	V(\tilde{\Th}, t) \\ & \ \ \ \ \quad \ge \frac{T}{2} \lambda_\mrm{min}(\bar{\Gamma}^{-1}) \|\tilde{\Th}(t)\|^2 > 0.
			\end{split}	
		\end{equation}   
		The time derivative of $V(\tilde{\Th}, t)$ along the solutions of \eqref{eq:err_dyn2} is obtained from
		\begin{equation}
			\label{eq:V_dot1}
			\begin{split}
				\dot{V}(\tilde{\Th}, t) 
				=\tilde{\Th}^\top(t)\bigg(&-\left(\bar{\Gamma}\la^\top(t)\la(t) \right)^\top P(t)\\& -P(t)\bar{\Gamma}\la^\top(t)\la(t)  +\dot{P}(t)\bigg)\tilde{\Th}(t).
			\end{split}
		\end{equation}
		Recall the following properties of the state-transition matrix of the system \eqref{eq:err_dyn2} (see, e.g., \cite[Theorem 1.1.1]{abou-kandil_matrix_2003})
		\begin{equation*}
			\begin{split}
				\bar{\Phi}(t,t) &= I_{nN}, \quad
				\frac{\partial \bar{\Phi}(s,t)}{\partial t} = \bar{\Phi}(s,t)\bar{\Gamma} \la^\top(t)\la(t), \\
				\frac{\partial \bar{\Phi}^\top(s,t)}{\partial t} &=  \la^\top(t)\la(t)\bar{\Gamma}\bar{\Phi}^\top(s,t),
			\end{split}
		\end{equation*}
		then $\dot{P}(t)$ can be expressed using the Leibniz integral rule as
		\begin{equation*}
			\begin{split}
				\dot{P}(t) &= \int_{t-T}^{t}\frac{\partial}{\partial t}(s-t+T)f(s,t) \mrm{d}s +  T \la^\top(t)\la(t) \\
				& = -\int_{t-T}^{t}f(s,t) \mrm{d}s +  T \la^\top(t)\la(t) \\& \ \ \ +  \la^\top(t)\la(t)\bar{\Gamma}\int_{t-T}^{t}(s-t+T)f(s,t) \mrm{d}s \\& \ \ \ +  \int_{t-T}^{t}(s-t+T)f(s,t)\mrm{d}s \bar{\Gamma} \la^\top(t)\la(t)
			\end{split}	
		\end{equation*}
		Consequently, $	\dot{V}(\tilde{\Th}, t)$ can be upper-bounded by substitution of $\dot{P}(t)$ in \eqref{eq:V_dot1} 
		\begin{equation}
			\label{eq:dot_V_bounds}
			\begin{split}
				&\dot{V}(\tilde{\Th}, t) =\tilde{\Th}^\top(t)\bigg(-\la^\top(t)\la(t)\bar{\Gamma} \int_{t-T}^{t}(s-t+T)f(s,t)\mathrm{d}s  
				\\& \quad - \int_{t-T}^{t}(s-t+T)f(s,t)\mathrm{d}s \bar{\Gamma}\la^\top(t)\la(t) -T\la^\top(t)\la(t)\\& \quad +  \la^\top(t)\la(t)\bar{\Gamma}\int_{t-T}^{t}(s-t+T)f(s,t) \mrm{d}s \\& \quad +  \int_{t-T}^{t}(s-t+T)f(s,t)\mrm{d}s \bar{\Gamma} \la^\top(t)\la(t)+  T \la^\top(t)\la(t)
				\\		&  \quad-\int_{t-T}^{t}f(s,t) \mrm{d}s \bigg)\tilde{\Th}(t)\\ & =- \tilde{\Th}^\top(t) \int_{t-T}^{t}f(s,t)\mrm{d}s\tilde{\Th}(t) \le -\Ml \|\tilde{\Th}(t)\|^2 <0.
			\end{split}
		\end{equation} 
		Hence, $V(\tilde{\Th}, t)$ defined in \eqref{eq:lya_fun} is a strong \ac{lf} for the system \eqref{eq:err_dyn2}. Utilizing Lyapunov's stability theorem (see, e.g., \cite[Theorem 4.1]{khalil_nonlinear_2002}), we can conclude \ac{guas} of the zero equilibrium solution of \eqref{eq:err_dyn2}.
	\end{proof}
	
	\subsection{Proof of Corollary \ref{corollary:convergence_bounds}}
	\begin{proof}[\unskip\nopunct]
		As seen from \eqref{eq:V_bounds} and \eqref{eq:dot_V_bounds}, the \ac{lf} given in Theorem~\ref{th:stong_lya} satisfies
		\begin{equation}
			\label{eq:bounds_lya}
			\begin{split}
				&	\kappa_1 \|\tilde{\Th}(t)\|^2 \ge	V(\tilde{\Th}, t) \ge \kappa_2  \|\tilde{\Th }(t)\|^2 , \\   
				&	\dot{V}(\tilde{\Th}, t) \le  -\Ml\|\tilde{\Th}(t)\|^2  ,
		\end{split}\end{equation} 
		with $$	\kappa_1 =	\left( \frac{T}{2}\lambda_\mrm{max}(\bar{\Gamma}^{-1}) + T\Mu\right), \quad \kappa_2 = \frac{T}{2} \lambda_\mrm{min}(\bar{\Gamma}^{-1})  .$$ 
		From \eqref{eq:bounds_lya} it follows that 
		\begin{equation*}
			\dot{V}(\tilde{\Th}, t) \le  -\frac{\Ml}{\kappa_1}V(\tilde{\Th},t),
		\end{equation*}
		and hence, from the Comparison Lemma \cite[Chap.~3]{khalil_nonlinear_2002}, we have		\begin{equation*}
			\int_{V(\tilde{\Th},t_0)}^{V(\tilde{\Th},t)} \frac{1}{v(\tilde{\Th},\tau)}	\mrm{d}v(\tilde{\Th},\tau) \le  	\int_{t_0}^{t}-\frac{\Ml}{\kappa_1}\mrm{d}\tau.
		\end{equation*}
		Solving the integrals provides
		\begin{equation*}
			{V}(\tilde{\Th}, t) \le  {V}(\tilde{\Th}, t_0)\mrm{e}^{-\frac{\Ml}{\kappa_1}(t-t_0)}.
		\end{equation*}
		Now, by substituting the right-hand side by the lower and the left-hand side by the upper bound provided in \eqref{eq:bounds_lya}, we get
		\begin{equation*}
			\kappa_2 \|\tilde{\Th}(t)\|^2 \le  	\kappa_1 \|\tilde{\Th}(t_0)\|^2\mrm{e}^{-\frac{ \Ml}{\kappa_1}(t-t_0)}.
		\end{equation*}
		The proof is completed by taking the square root and rearranging the terms to obtain 
		$$ 	 \|\tilde{\Th}(t)\| \le  	\sqrt{\frac{\kappa_1}{\kappa_2}} \|\tilde{\Th}(t_0)\|\mrm{e}^{-\frac{ \Ml}{2\kappa_1}(t-t_0)}.$$
	\end{proof}
	
	\subsection{Proof of Corollary \ref{corollary:ISS_bounds}}
	\begin{proof}[\unskip\nopunct]
		The time derivative of $V(\tilde{\Th}, t)$ \eqref{eq:lya_fun} along the solutions of \eqref{eq:disturbed_system_ISS} is obtained from
		\begin{equation*}
			\begin{split}
				\dot{V}(\tilde{\Th}, t) 
				&=\tilde{\Th}^\top(t)\bigg(-\la^\top(t)\la(t)\bar{\Gamma}P(t) -P(t)\bar{\Gamma}\la^\top(t)\la(t) \\& +\dot{P}(t)\bigg)\tilde{\Th}(t) + \tilde{\Th}^\top(t)  P(t)  \delta_\mrm{ISS}(t) + \delta_\mrm{ISS}^\top(t)  P(t) \tilde{\Th}(t)   
				\\	& \le  -\Ml\|\tilde{\Th}(t)\|^2 + 2 \kappa_1 \|\tilde{\Th}(t)\| \|\delta_\mrm{ISS}(t)\|. 
			\end{split}
		\end{equation*}
		We use the term $ -\Ml\|\tilde{\Th}(t)\|^2$ to dominate $2 \kappa_1 \|\tilde{\Th}(t)\| \|\delta_\mrm{ISS}(t)\|$ for sufficiently large $\|\tilde{\Th}(t)\|$
		\begin{equation*}
			\begin{split}
				\dot{V}(\tilde{\Th}, t) 
				\le 	& -\Ml(1-\beta)\|\tilde{\Th}(t)\|^2 + \bigg( -\Ml\beta\|\tilde{\Th}(t)\| \\ &+ 2 \kappa_1  \|\delta_\mrm{ISS}(t)\| \bigg) \|\tilde{\Th}(t)\|,
			\end{split}
		\end{equation*}
		where we introduced $0<\beta<1$. Hence, we have 
		\begin{equation*}
			\dot{V}(\tilde{\Th}, t) 
			\le  -\Ml(1-\beta)\|\tilde{\Th}(t)\|^2, \ \forall  \  \|\tilde{\Th}(t)\| \ge \frac{2 \kappa_1 }{\Ml \beta}\| \|\delta_\mrm{ISS}(t)\| ,
		\end{equation*}
		fulfilling the conditions of \cite[Theorem 4.19]{khalil_nonlinear_2002} with \begin{equation*}
			\begin{split}
				&\alpha_1\left(\|\tilde{\Th}(t)\|^2\right) = \kappa_2 \|\tilde{\Th}(t)\|^2, \ \ \ \alpha_2\left(\|\tilde{\Th}(t)\|^2\right) = \kappa_1 \|\tilde{\Th}(t)\|^2, \\ &\rho\left(\|\delta_\mrm{ISS}(t)\|\right) = \frac{2 \kappa_1 }{\Ml\beta}\| \|\delta_\mrm{ISS}(t)\|.
			\end{split}
		\end{equation*}
		According to \cite[Theorem 4.19]{khalil_nonlinear_2002}, the \ac{iss}-gain from the disturbance $\delta_\mrm{ISS}(t)$ to the state $\tilde{\Th}(t)$ is
		\begin{equation*}
			\begin{split}
				\gamma_{\mrm{ISS}}\left(\|\delta_\mrm{ISS}(t)\|\right) &= \alpha_1^{-1} \circ \alpha_2 \circ \rho\left(\|\delta_\mrm{ISS}(t)\|\right) \\&=  \frac{2\kappa_1}{\Ml\beta}\sqrt{\frac{\kappa_1}{\kappa_2}}\left(\|\delta_\mrm{ISS}(t)\|\right), 
			\end{split}
		\end{equation*} 
		where $\circ$ denotes the function composition operator. Hence, the upper bound of the ISS gain is obtained as denoted in \eqref{eq:ISS_gain}. 
	\end{proof}

	\subsection{Proof of Theorem \ref{th:tuning}}
	\begin{proof}[\unskip\nopunct]
		We modify the strong \ac{lf} from Theorem~\ref{th:stong_lya}, by introducing the scalar weighing factors $c_1>0$ and $c_2>0$
		\begin{equation*}
			V^\mrm{d}(\tilde{\Th}^\mrm{d}, t) = (\tilde{\Th}^\mrm{d})^\top(t)P^\mrm{d}(t)\tilde{\Th}^\mrm{d}(t),
		\end{equation*}
		with 
		\begin{equation*}
			\begin{split}
				P^\mrm{d}(t) &= c_1\bar{\Gamma}^{-1} + c_2\bar{\Pi}(t), \\ \bar{\Pi}(t) & = \int_{t-T}^{t}(s-t+T)\bar{\Phi}^\top(s,t)\la^\top(s)\la(s)\bar{\Phi}(s,t)\mathrm{d}s.
			\end{split}
		\end{equation*}
		The time derivative of $V^\mrm{d}(\tilde{\Th}^\mrm{d}, t)$ along the solutions of \eqref{eq:err_dyn_dis} is obtained from\footnote{For the remainder of the notation, all signals' time arguments are omitted.}
		\begin{equation*}
			\begin{split}
				&	\dot{V}^\mrm{d}(\tilde{\Th}^\mrm{d}, t) 
				=  (\tilde{\Th}^\mrm{d})^\top\left( \left(c_2T-2c_1 \right) \la^\top  \la -c_2 \bar{\Pi}\right) \tilde{\Th}^\mrm{d} \\&
				+  (\tilde{\Th}^\mrm{d})^\top\left(c_1 \la^\top\Db2 +c_2\bar{\Pi}\bar{\Gamma} \la^\top\Db2 -c_1 \bar{\Gamma}^{-1}\Db1 -c_2 \bar{\Pi}\Db1\right) \delta
				\\&
				+  \delta^\top\left(c_1 \Db2^\top\la +c_2 \Db2^\top\la\bar{\Gamma} \bar{\Pi}-c_1 \Db1^\top\bar{\Gamma}^{-1} -c_2 \Db1^\top\bar{\Pi}\right) \tilde{\Th}^\mrm{d}.
				%
				%
				%
			\end{split}
		\end{equation*}
		We get the following set of inequalities by using Young's inequality
		\begin{equation*}
			\begin{split}
				&c_1(\tilde{\Th}^\mrm{d})^\top\la^\top \Db2 \delta + c_1\delta^\top \Db2^\top \la \tilde{\Th}^\mrm{d}    \le  \\  & \qquad \qquad \qquad c_1 \tilde{\Th}^\mrm{d}\la^\top \la \tilde{\Th}^\mrm{d}+ c_1\delta^\top \Db2^\top \Db2 \delta \ , \\
				&	c_1	(\tilde{\Th}^\mrm{d})^\top\bar{\Gamma}^{-1}\Db1{\delta} +c_1 \delta^\top \Db1^\top\bar{\Gamma}^{-1}{\tilde{\Th}^\mrm{d}}   \le \\ & \qquad \qquad \qquad \frac{c_1}{\epsilon_1}\delta^\top \Db1^\top\Db1 \delta +c_1\epsilon_1(\tilde{\Th}^\mrm{d})^\top \left(\bar{\Gamma}^{-1}\right)^2 \tilde{\Th}^\mrm{d} \ , \\
				&	c_2	(\tilde{\Th}^\mrm{d})^\top\bar{\Pi}\bar{\Gamma}\la^\top \Db2{\delta} +c_2 \delta^\top \Db2^\top\la\bar{\Gamma}\bar{\Pi}{\tilde{\Th}^\mrm{d}}   \le \\ & \qquad \qquad \qquad\frac{c_2}{\epsilon_2}(\tilde{\Th}^\mrm{d})^\top\bar{\Pi} \bar{\Gamma}\la^{\top} \la \bar{\Gamma}\bar{\Pi}\tilde{\Th}^\mrm{d} +c_2\epsilon_2\delta^\top \Db2^\top\Db2\delta \ , \\ 
				&	c_2	(\tilde{\Th}^\mrm{d})^\top\bar{\Pi}\Db1{\delta} +c_2 \delta^\top \Db1\bar{\Pi}{\tilde{\Th}^\mrm{d}} \le\\&  \qquad \qquad \qquad \frac{c_2}{\epsilon_3}\delta^\top \Db1^\top\Db1 \delta +c_2\epsilon_3(\tilde{\Th}^\mrm{d})^\top \bar{\Pi}^2\tilde{\Th}^\mrm{d} .
			\end{split}
		\end{equation*}
		Defining $\chi^\top=\begin{bmatrix}(\tilde{\Th}^\mrm{d})^\top \la^\top &(\tilde{\Th}^\mrm{d})^\top & z^\top & \delta^\top
		\end{bmatrix}$ and 
		\begin{equation*}
			\begin{split}
				\epsilon_1 &= \frac{\lambda^2_\mrm{min}(\bar{\Gamma})c_2\Ml}{4c_1} \ , \quad  \epsilon_2 = \frac{8\Mu^2r_4\lambda^2_\mrm{max}(\bar{\Gamma})}{\Ml} \ , \quad 	\epsilon_3= \frac{\Ml}{8\Mu^2},
			\end{split}
		\end{equation*}
		we bound the time derivative of $V^\mrm{d}(\tilde{\Th}^\mrm{d})$ as
		\begin{equation}
			\label{eq:Vd_dot}
			\dot{V}^\mrm{d}(\tilde{\Th}^\mrm{d}) \le \chi^\top \Psi_1\chi + \frac{4c_1^2}{c_2\lambda^2_\mrm{min}(\bar{\Gamma})\Ml} \delta^\top \Db1^\top\Db1 \delta,
		\end{equation}
		with $\Psi_1$ defined in \eqref{eq:Psi}.
		\begin{figure*}[!t]
			\begin{equation}
				\scriptsize
				\label{eq:Psi}
				\begin{aligned}
					\Psi_{1}&= 	\begin{bmatrix}\psi_{11}  & 0 & 0 & 0 \\ 
						0 & \frac{-c_2\Ml}{2}I_{nN} & 0 & 0 \\
						0 & 0 & 0 & 0 \\
						0 & 0 & 0  &  \psi_{44}  \end{bmatrix} \ , \quad 
					\Psi_{2}= 	\begin{bmatrix}0  & 0 & 0 & 0 \\ 
						0 & 0 & Q^\top & 0 \\
						0 & Q & -I_{p} & W \\
						0 & 0 & W^\top  &  -\gamma I_{r}  \end{bmatrix} \ , \quad
					\Psi_{2}^\mrm{OE}= 	\begin{bmatrix}0  & 0 & Q^\top & 0 \\ 
						0 & 0 & 0 & 0 \\
						Q & 0 & -I_{p} & W^\mrm{OE} \\
						0 & 0 & (W^\mrm{OE})^\top  &  -\gamma I_{r}  \end{bmatrix} \ , \\
					\psi_{11} &=    \left(-c_1 + c_2T \right)I_{\left(nN_\mrm{y}+Nn_{\mrm{e}(t)}\right)}   \ ,  \quad \psi_{44} = \frac{8c_2\Mu^2}{\Ml} \Db1^\top \Db1 +	 \left(c_1 +  \frac{8c_2\Mu^2r_4\lambda^2_\mrm{max}(\bar{\Gamma})}{\Ml}\right) \Db2^\top \Db2     .
				\end{aligned}
			\end{equation}
		\hrule 
		\end{figure*}
		To investigate the $\mathcal{L}_2$-gain of the error system \eqref{eq:err_dyn_dis}, we recall that if $\dot{V}^\mrm{d}(\tilde{\Th}^\mrm{d}) \le -z^\top z +\gamma \delta^\top \delta$ holds the system \eqref{eq:err_dyn_dis} has an $\mathcal{L}_2$-gain of less or equal to $\sqrt{\gamma}$ \cite[Section 8.1]{van_der_schaft_l2-gain_2017}. Using \eqref{eq:Vd_dot} and \eqref{eq:zTz}, we can rewrite this condition as
		\begin{equation}
			\label{eq:Vd_dot2}
			\begin{split}
				\dot{V}^\mrm{d}(\tilde{\Th}^\mrm{d}) &+z^\top z -\gamma \delta^\top \delta \le \chi^\top \Psi_1\chi + \frac{4c_1^2}{c_2\lambda^2_\mrm{min}(\bar{\Gamma})\Ml} \delta^\top \Db1^\top\Db1 \delta \\
				& \qquad \quad \ \qquad \qquad+z^\top z  -\gamma \delta^\top \delta +z^\top z - z^\top z
				\\ & = \chi^\top \left(\Psi_1+\Psi_2  +\Phi_{12}\left(\lambda^2_\mrm{min}(\bar{\Gamma}) I_{nN} \right)^{-1} \Phi_{12}^\top \right) \chi \le0,
			\end{split}
		\end{equation}
		with $\Psi_2$ and $\Phi_{12}$ defined in \eqref{eq:Psi} and \eqref{eq:Phi}, respectively. Utilizing the Schur complement \cite[Sec. 4.2.3]{herrmann_linear_2007-1}, we can verify that condition \eqref{eq:Vd_dot2} is fulfilled for 
		\begin{equation*}
			\begin{bmatrix}
				\Psi_{1}+\Psi_2 & \Phi_{12}\\ \Phi_{12}^\top & -\lambda^2_\mrm{min}(\bar{\Gamma}) I_{nN}\end{bmatrix} \le0.
		\end{equation*}
		Defining $\gamma_1 \coloneqq \lambda^2_\mrm{max}(\bar{\Gamma})$ and  $\gamma_2 \coloneqq \lambda^2_\mrm{min}(\bar{\Gamma})$ and noticing that $\Phi_{11} = \Psi_{1}+\Psi_2$ completes the proof.
	\end{proof}
	
	\subsection{Proof of Corollary \ref{co:tuning2}}
	\begin{proof}[\unskip\nopunct]
		We 	define $$(\chi^\mrm{OE})^\top=\begin{bmatrix}(\tilde{\Th}^\mrm{d})^\top \la^\top &(\tilde{\Th}^\mrm{d})^\top & (z^\mrm{OE})^\top & \delta^\top
		\end{bmatrix},$$ and follow the proof of Theorem~\ref{th:tuning} until \eqref{eq:Vd_dot2}, where we replace $z$ with the modified performance output $z^\mrm{OE}$ defined in \eqref{eq:z_OE} and $\chi$ with $\chi^\mrm{OE}$. This yields 
		\begin{equation}
			\label{eq:Vd_dot3}
			\begin{split}
				&\dot{V}^\mrm{d}(\tilde{\Th}^\mrm{d}) +(z^\mrm{OE})^\top z^\mrm{OE} -\gamma \delta^\top \delta \le  (\chi^\mrm{OE})^\top \Psi_1\chi^\mrm{OE} \\& \quad+ \frac{4c_1^2}{c_2\gamma_2\Ml} \delta^\top \Db1^\top\Db1 \delta +(z^\mrm{OE})^\top z^\mrm{OE}  -\gamma \delta^\top \delta \\&\quad+(z^\mrm{OE})^\top z^\mrm{OE} - (z^\mrm{OE})^\top z^\mrm{OE}
				\\ & = (\chi^\mrm{OE})^\top \left(\Psi_1+\Psi_2^\mrm{OE}  +\Phi_{12}\left(\lambda^2_\mrm{min}(\bar{\Gamma}) I_{nN} \right)^{-1} \Phi_{12}^\top \right) \chi^\mrm{OE} \\&\le0,
			\end{split}
		\end{equation}
		with $\Psi_1$ and $\Psi_2^\mrm{OE}$ defined in \eqref{eq:Psi} and $\Phi_{12}$ defined in \eqref{eq:Phi}.
		Utilizing the Schur complement \cite[Sec. 4.2.3]{herrmann_linear_2007-1}, we can verify that condition \eqref{eq:Vd_dot3} is fulfilled for 
		\begin{equation*}
			\begin{bmatrix}
				\Psi_{1}+\Psi_2^\mrm{OE} & \Phi_{12}\\ \Phi_{12}^\top & -\gamma_2 I_{nN}\end{bmatrix} \le0.
		\end{equation*}
		Noticing that $\Phi_{11}^\mrm{OE} = \Psi_{1}+\Psi_2^\mrm{OE}$ completes the proof.
	\end{proof}
	
	\subsection{Proof of Proposition \ref{pro:tuning}}
	\begin{proof}[\unskip\nopunct]
		From \eqref{eq:Vd_dot2}, we define 
		\begin{equation*}
			\begin{split}
				\Omega(\gamma, \gamma_1, \gamma_2, &\alpha, c_2, c_1, \Mu, \Ml,r_4 ,T)\coloneqq \\& \chi^\top \left(\Phi_{11}  +\Phi_{12}\left(\gamma_2 I_{nN} \right)^{-1} \Phi_{12}^\top \right) \chi ,
			\end{split}
		\end{equation*}
		and 		$\Omega^1 \coloneqq\Omega(\gamma^*, \gamma_1^*, \gamma_2^*, \alpha^*, c_1^*, c_2^*, \Mu, \Ml,r_4 ,T)$. Replacing the upper and lower bounds of the excitation parameters \eqref{eq:bounds_excitation}, we have $\Omega^* \coloneqq\Omega(\gamma^*, \gamma_1^*, \gamma_2^*, \alpha^*, c_1^*, c_2^*, \Mu^*, \Ml^*,r_4^* ,T^*) \le 0,$
		since \eqref{eq:LMI_tuning_th1} is fulfilled as defined in Proposition~\ref{pro:tuning} and \eqref{eq:LMI_tuning_th1} implies \eqref{eq:Vd_dot2}.  
		\begin{figure*}[!t]
			\begin{equation}
				\scriptsize
				\label{eq:tilde_Phi}
				\begin{aligned}
					\tilde{\Phi}_{11}&= 	\begin{bmatrix} c_2^*(T-T^*) I_{\left(nN_\mrm{y}+Nn_{\mrm{e}(t)}\right)}  & 0 & 0 & 0 \\ 
						0 & \frac{-c_2^*(\Ml-\Ml^*)}{2}I_{nN} & 0 & 0 \\
						0 & 0 & 0 & 0 \\
						0 & 0 & 0  &  8c_2^*\left(\frac{\Mu^2}{\Ml}-\frac{(\Mu^*)^2}{\Ml^*}\right) \Db1^\top \Db1 +	 8c_2^*\gamma_1^*\left(\frac{\Mu^2r_4}{\Ml}- \frac{(\Mu^*)^2r_4^*}{\Ml^*}\right) \Db2^\top \Db2      \end{bmatrix} \ . 
				\end{aligned}
			\end{equation}
			\hrule 
		\end{figure*}
		We evaluate the difference 
		\begin{equation*}
			\begin{split}
				\Omega^1-\Omega^* =\chi^\top \left(\tilde{\Phi}_{11}  \right) \chi  +\frac{4(c_1^*)^2}{c_2^*\gamma_2^*}\left(\frac{1}{\Ml}-\frac{1}{\Ml^*}\right) \delta^\top \Db1^\top \Db1 \delta,
			\end{split}
		\end{equation*}
		where $\tilde{\Phi}_{11}$ is defined in \eqref{eq:tilde_Phi}. From \eqref{eq:bounds_excitation}, we have $$ T^*\ge T, \ \frac{1}{\Ml^*} \ge \frac{1}{\Ml} , \ \frac{(\Mu^*)^2r_4^*}{\Ml^*} \ge \frac{\Mu^2r_4}{\Ml}, \   \frac{(\Mu^*)^2}{\Ml^*} \ge   \frac{\Mu^2}{\Ml}, $$
		and, hence, $\Omega^1-\Omega^* \le 0$. Consequently, we conclude that
		\begin{equation*}
			\dot{V}^\mrm{d}(\tilde{\Th}^\mrm{d}) +z^\top z -\gamma \delta^\top \delta \le \Omega^1 \le 0,
		\end{equation*}
		which implies \eqref{eq:LMI_tuning_th1} using the Schur complement.
	\end{proof}
	\section*{Acknowledgements} This research is partly supported by the German Federal Government, the Federal Ministry of Education and Research and the State of Brandenburg within the framework of the joint project EIZ: Energy Innovation Center (project numbers 85056897 and 03SF0693A) with funds from the Structural Development Act (Strukturstärkungsgesetz) for coal-mining regions.
	\bibliographystyle{IEEEtran}
	\bibliography{Diss_bib}
	
	\begin{IEEEbiography}[{\includegraphics[width=1in,height=1.25in,clip,keepaspectratio]{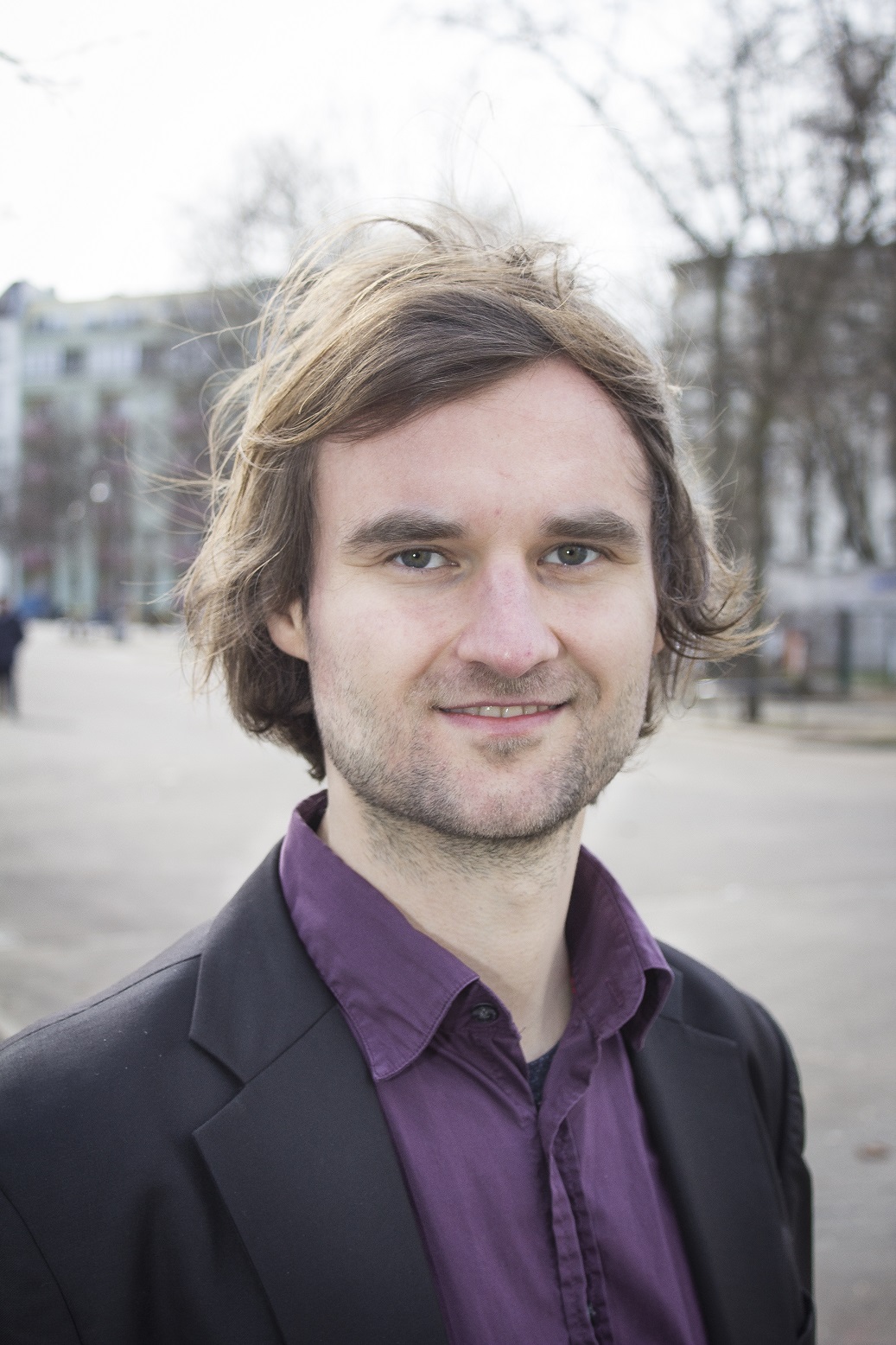}}]{Nicolai Lorenz-Meyer}
		received his M.Sc. in
		Engineering Science from the Technical University
		of Berlin in 2019 and his B.Eng. in Business Administration and Engineering for Environment and
		Sustainability from the Beuth University of Applied
		Sciences and the Berlin School of Economics and
		Law in 2016.
		He is currently pursuing the Ph.D. degree with the
		the chair of Control Systems and Network Control
		Technology at the Brandenburg University of Technology Cottbus-Senftenberg, Germany. His current research interests include the development of control theory-based methods, especially in the context of distributed parameter estimation, and their application to online dynamics security assessment in power systems.
	\end{IEEEbiography}
	\begin{IEEEbiography}[{\includegraphics[width=1in,height=1.25in,clip,keepaspectratio]{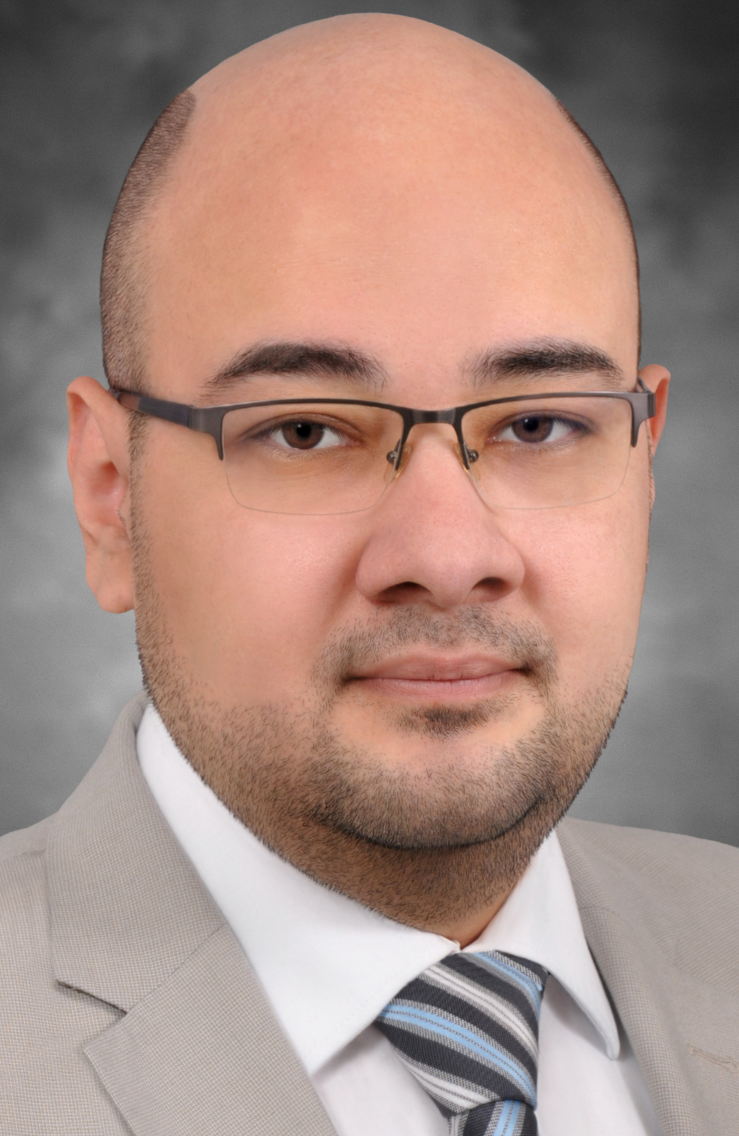}}]
		{Juan G. Rueda-Escobedo} received his Licentiate in Mechatronic Engineering and both his Master's and Doctoral degrees (with honors) in Electrical Engineering, specializing in Automatic Control, from the National Autonomous University of Mexico (UNAM), Mexico City, in 2013, 2014, and 2018, respectively. He served as a Research Associate at the Chair of Control Systems and Network Control Technology at Brandenburg University of Technology Cottbus-Senftenberg, Germany, from 2018 to 2023. Since 2023, he has been appointed as an Associate Professor in the Basic Sciences Division of the Engineering Faculty at UNAM. His research primarily focuses on the design of observers for time-varying and nonlinear systems, development of innovative parameter identification strategies, and the advancement of adaptive observation and control techniques with enhanced convergence rates.
	\end{IEEEbiography}
	\begin{IEEEbiography}[{\includegraphics[width=1in,height=1.25in,clip,keepaspectratio]{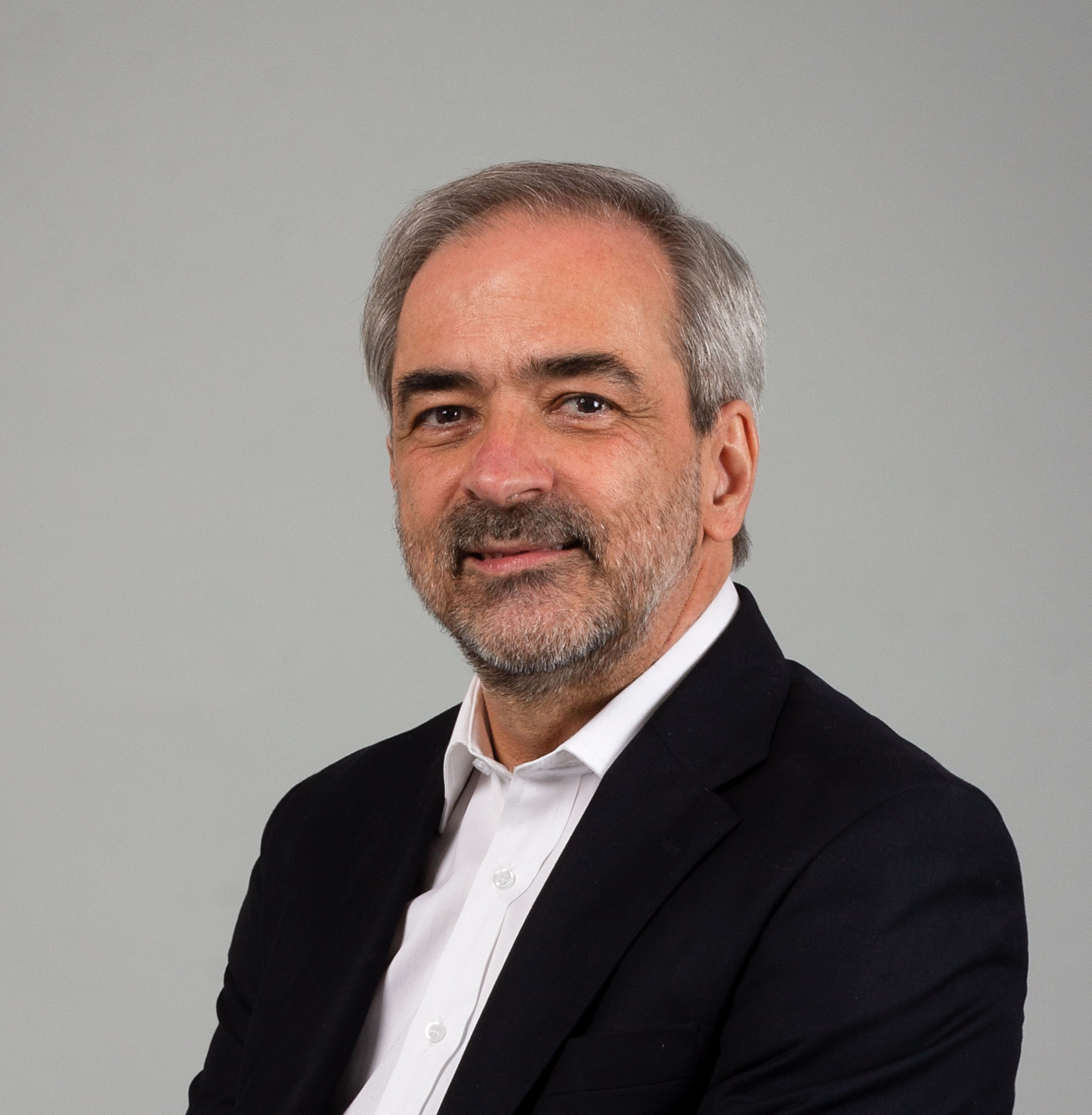}}]
		{Jaime~A.~Moreno} was born in Colombia and he received his PhD degree (Summa cum Laude) in Electrical Engineering (Automatic Control) from the Helmut-Schmidt University in Hamburg, Germany in 1995. The Diplom-Degree in Electrical Engineering (Automatic Control) from the Universität zu Karlsruhe, Karlsruhe, Germany in 1990, and the Licentiate-Degree (with honors) in Electronic Engineering from the Universidad Pontificia Bolivariana, Medellin, Colombia in 1987. He is full Professor of Automatic Control and he was the Head of the Electrical and Computing Department at the Institute of Engineering from the National University of Mexico (UNAM), in Mexico City. He has served IFAC in different capacities. He is the author and editor of 8 books, 14 book chapters, 1 patent, and author and co-author of more than 500 papers in refereed journals and conference proceedings. His current research interests include robust and non-linear control, in particular, with emphasis on Lyapunov methods for higher order sliding modes control, with applications to biochemical (wastewater treatment processes) and electromechanical processes, and the design of nonlinear observers.
	\end{IEEEbiography}
	
	\begin{IEEEbiography}[{\includegraphics[width=1in,height=1.25in,clip,keepaspectratio]{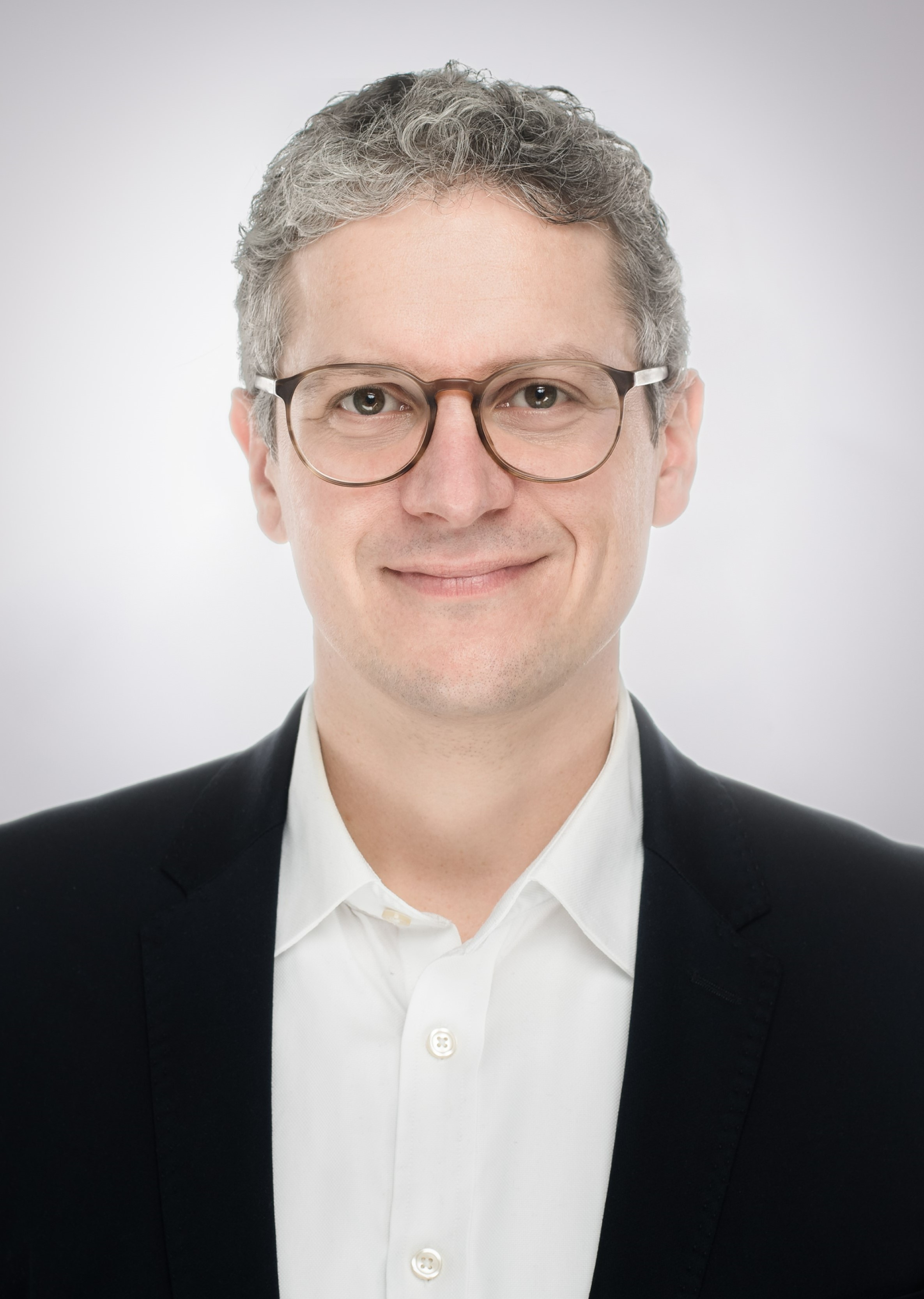}}]
		{Johannes Schiffer} received the Diploma degree in engineering cybernetics from the University of Stuttgart, Germany, in 2009 and the Ph.D. degree (Dr.-Ing.) in electrical engineering from Technische Universität (TU) Berlin, Germany, in 2015. He currently holds the chair of Control Systems and Network Control Technology at Brandenburgische Technische Universität Cottbus-Senftenberg, Germany and leads the business area Control, Automation and Operation Management at the Fraunhofer Research Institution for Energy Infrastructures and Geothermal Systems (IEG). Prior to that, he has held appointments as Lecturer (Assistant Professor) at the School of Electronic and Electrical Engineering, University of Leeds, U.K. and as Research Associate in the Control Systems Group and at the Chair of Sustainable Electric Networks and Sources of Energy both at TU Berlin. He and his co-workers received the Automatica Paper Prize over the years 2014 - 2016 and the at - Automatisierungstechnik Best Paper Award 2022. His current research interests include the development of control engineering methods for complex and interconnected dynamic systems and their application to practical challenges, in particular to increase climate neutrality in industry and society.
	\end{IEEEbiography}
\end{document}